\documentclass[proc]{edpsmath}
\usepackage{graphicx} 
\usepackage{fullpage}

\usepackage{comment}
\usepackage{xcolor}
\usepackage{url}

\usepackage{amsmath}
\usepackage{amsthm}
\usepackage{amssymb}
\usepackage{stmaryrd}

\usepackage{graphicx}
\usepackage{algorithm}
\usepackage{algorithmic}

\usepackage{verbatim}
\usepackage{pgf,tikz}
\usepackage{subcaption}
\usepackage{tcolorbox}
\usepackage{hyperref}

\usepackage{setspace}

\newcommand{\norm}[1]{\left \| #1 \right \|}
\newcommand{\abs}[1]{\left | #1 \right |}
\newcommand{\set}[2]{\left \{ #1 \, : \, #2 \right \}}

\DeclareMathOperator{\R}{\mathbb{R}} 
\DeclareMathOperator{\Z}{\mathbb{Z}} 
\DeclareMathOperator{\N}{\mathbb{N}}

\newcommand{\bfm}[1]{\mathbf{#1}}

\newcommand{\red}[1]{\textcolor{red}{#1}}

\newtheorem{proposition}{Proposition}
\newtheorem{remark}{Remark}
\begin{document}

\title{Comparison between tensor methods and 
neural networks in electronic structure calculations.
}\thanks{The authors acknowledge funding from the Tremplin-ERC Starting ANR grant HighLEAP (ANR-
22-ERCS-0012). This publication is part of a project that has received funding from the European
Research Council (ERC) under the European Union’s Horizon 2020 Research and Innovation Programme – Grant Agreement n ◦ 101077204.}\thanks{This work
was supported by the French ‘Investissements d’Avenir’ program, project Agence
Nationale de la Recherche (ISITE-BFC) (contract ANR-15-IDEX-0003). GD was also supported by the Ecole des Ponts-ParisTech. 
GD aknowledges the support of the region Bourgogne Franche-Comt\'e.}
\renewcommand{\authors}[3][]{
      \stepcounter{ifmbe@authors}
      \expandafter\def\csname ifmbe@author\alph{ifmbe@authors}\endcsname
      {#2$^{\expandafter\the\csname ifmbe@affiliationcounter#3\endcsname
        \if\relax\detokenize{#1}\relax\else,#1\fi}$}
}

\author{Mathias Dus}\address{École nationale des ponts et chaussées and Inria Paris team Matherials, Paris, France}
\author{Geneviève Dusson}\address{Laboratoire de Mathématiques, Université de Franche-Comté, Besançon, France}
\author{Virginie Ehrlacher}\address{École nationale des ponts et chaussées and Inria Paris team Matherials, Paris, France}
\author{Clément Guillot}\address{École nationale des ponts et chaussées and Inria Paris team Matherials, Paris, France}
\author{Joel Pascal Soffo Wambo$^{5, }$}
\address{Centrale Nantes, Nantes Université, Laboratoire de Mathématiques Jean Leray, Nantes, France}
\address{Airbus Central Research \& Technology, Virtual Product Engineering/XRV}
\begin{abstract} 
This article compares the tensor method density matrix renormalization group (DMRG) with two neural network based methods -namely FermiNet and PauliNet) for determining the ground state wavefunction of the  many-body electronic Schrödinger problem. We provide numerical simulations illustrating the main features of the methods and showing convergence with respect to some parameter, such as the rank for DMRG, and number of pretraining iterations for neural networks. We then compare the obtained energy with the methods for a few atoms and molecules, for some of which the exact value of the energy is known for the sake of comparison. In the last part of the article, we propose a new kind of neural network to solve the Schr\"odinger problem based on the training of the wavefunction on a simplex, and an explicit permutation for evaluating the wavefunction on the whole space. We provide numerical results on a toy problem for the sake of illustration.
\end{abstract}

%

\maketitle


\section*{Introduction}

Electronic structure calculations are used in many scientific fields, in order to accurately compute physical properties of molecular of materials systems. In this field, a common task consists in computing a numerical approximation of the state of the set of electrons in a molecule of interest with lowest possible energy. This is called the \itshape ground state \normalfont of the electrons. This is of vital importance in many applications since the knowledge of the ground state enbales to obtain predictions on many physical properties of the molecule, like electric, chemical, optical or magnetic properties. 

More precisely, a molecular system is in general composed of $M$ nuclei (composed of protons and netrons) and $N$ electrons. Due to the large mass of the nuclei compared to the electrons, the latter are often modeled as classical particles, while electrons are considered as quantum particles. This is called the Born--Oppenheimer approximation~\cite{Born1927-ap}. Therefore, the nuclei are described by their electric charges and their positions and velocities in space, while the electrons are described by a wavefunction $\psi: (\mathbb{R}^3 \times \mathbb{Z}^2)^N \to \mathbb{C}$, which depends on all positions in $\R^3$ and spins in $\Z_2$ of the $N$ electrons.  Since electrons are fermionic paricles, due to the Pauli exclusion principle, the wavefunction is antisymmetric with respect to permutation of the ordering of the electrons, which means that, for almost all $\bfm{r}_1, \ldots, \bfm{r}_N\in \mathbb{R}^3$ and all $\sigma_1, \ldots, \sigma_N \in \mathbb{Z}^2$,
\[
\forall s \in S_N, \quad \psi(\bfm{r}_{s(1)}, \sigma_{s(1)}, ..., \bfm{r}_{s(N)}, \sigma_{s(N)}) = (-1)^{\varepsilon(s)} \psi(\bfm{r}_1, \sigma_1, ..., \bfm{r}_N, \sigma_N),
\]
where $S_N$ is the set of permutations of the set $\{1, \ldots, N\}$. Any admissible wavefunction describing the state of the set of $N$ electrons has therefore to be an element of the space 
\[
\bigwedge_{n=1}^N L^2(\mathbb R^3 \times \mathbb Z_2),
\]
which is the anti-symmetric tensor product of $L^2(\mathbb R^3 \times \mathbb Z_2)$. From a physical point of view, the quantity $|\psi(\bfm{r}_{s(1)}, \sigma_{s(1)}, ..., \bfm{r}_{s(N)}, \sigma_{s(N)})|^2$ represents the probability density of finding, for any $1\leq i \leq N$, the $i^{th}$ electron at position $\bfm{r}_i\in \mathbb{R}^3$ with spin $\sigma_i \in \mathbb{Z}_2$. This implies in particular that the wavefunction $\psi$ also has to satisfy the following normalization condition: 
$$
\left\|\psi\right\|_{\bigwedge_{n=1}^N L^2(\mathbb R^3 \times \mathbb Z_2)} = 1. 
$$

As mentioned above, one key problem in the field is to find the electronic ground state of the molecular system, which minimizes the energy of the system. It requires to solve the so-called electronic many-body Schr\"odinger equation, an eigenvalue partial differential equation in very high dimension. Several methods are used to compute approximations of the ground state, including traditional wave-function methods, such as Hartree--Fock, Configuration Interaction, or Coupled Cluster, as well as Density Functional Theory (DFT), (see~\cite{Helgaker2014-py} for an introduction to these methods) or even the more recently developed methods based on tensor decomposition~\cite{Legeza_2014} and neural networks~\cite{schnet,Hermann_2020,Spencer2021-ee}. 

The aim of this article is to compare several methods which are claimed to be very accurate and efficient, but are, to our knowledge, never been compared in a common framework. Therefore, in this article, we provide numerical simulations based on the Density Matrix Renormalization Group (DMRG) method with neural network based methods, in particular the FermiNet~\cite{Spencer2021-ee} and PauliNet~\cite{Hermann_2020} approach. 

On top of this comparison, we propose a new type of neural network method to compute the ground state based on the learning of the wavefunction on the principal simplex of the spatial domain and explicit antisymmetrization for the computation of the wavefunction on the rest of the domain that we call hereafter SimplexNet.

The outline of this article is as follows. In section~\ref{sec:sec_approx}, we present the Schr\"odinger problem together with a few approximations of interest. In section~\ref{sec:num}, we provide numerical results comparing these methods. In section~\ref{sec:simplexnet}, we present the SimplexNet.
We finally provide a conclusion in section~\ref{sec:concl}.

\section{Many-body Schr\"odinger equation and some approximations}
\label{sec:sec_approx}

As mentioned in the introduction, we are interested in the solution of the ground state many-body Schr\"odinger equation. We present in this section the Schr\"odinger equation as well as three major approximations: the Hartree--Fock method, the Density Matrix Renormalization Group (DMRG) method, and Neural Networks based methods including FermiNet and PauliNet.

For the sake of simplicity, we choose to omit the spin variables in the following presentation, assuming that they are fixed. However, numerical simulations presented in Section~\ref{sec: Numerical experiments} do include spin variables.

\subsection{Many-body electronic Schr\"odinger equation}

Let us start by presenting the many-body electronic Schr\"odinger problem.  Given a molecular system with $M$ nuclei having positions $(\bfm{X}_m)_{1\le m \le M} \in \R^{3M}$ and charges $(Z_m)_{1\le m \le M} \in \N^M$, we consider the following electronic Hamiltonian given by
    \begin{equation}
        H
        =   -\frac 1 2 \sum_{n=1}^N \Delta_{\bfm{r}_n} + \sum_{n=1}^N \sum_{m=1}^M \frac{-Z_m}{\abs{\bfm{r}_n - \bfm{X}_m}} + \sum_{1 \leq j < n \leq N} \frac{1}{\abs{\bfm{r}_j - \bfm{r}_n}}.
    \end{equation}
    Under appropriate assumptions which we assumed to be satisfied here, from Zhislin's theorem, it holds that the Hamiltonian $H$ is a self-adjoint bounded from below operator acting upon the space of antisymmetric wave functions $\bigwedge_{j=1}^N L^2(\mathbb R^3)$ with domain $D:=H^2((\mathbb{R}^3)^N) \cap \bigwedge_{j=1}^N L^2(\mathbb R^3)$ and form domain $\Sigma:=H^1((\mathbb{R}^3)^N) 
 \cap \bigwedge_{j=1}^N L^2(\mathbb R^3)$. Denoting by $\langle \cdot, \cdot\rangle$ (and by $\|\cdot\|$) the canonical $L^2((\mathbb{R}^3)^N)$ scalar product (and its associated norm respectively), the electronic many-body Schr\"odinger equation then consists in solving the following minimization problem
        \begin{equation} \label{equ: minmax principle}
        E_0
        :=   \inf_{\psi \in \Sigma, \norm \psi = 1} \langle \psi, H \psi \rangle.
    \end{equation}
   Writing the Euler--Lagrange equations of the problem, the minimization problem can be written as the following eigenvalue problem: find $(E_0, \psi_0)\in \mathbb{R}\times \Sigma$ with $\|\psi_0\| =1$ such that
    \begin{equation}
    \label{eq:schrod}
             H\psi_0 = E_0 \psi_0.
    \end{equation}
      The value $E_0$ is often called the ground state energy of $H$, and any associated eigenvector $\psi_0$ is called a ground state.

    In the case of the hydrogen atom, where $N=M=1$ and $Z_1 = 1$, it is well known that the essential spectrum $\sigma_{\rm ess}(H)$ of $H$ is equal to $[0, \infty)$, and the discrete spectrum $\sigma_{\rm disc}(H) = \sigma(H) \cap (-\infty, 0)$ consists of an increasing sequence $(E_n)_{n \geq 0}$ with $\lim_{n \to \infty} E_n = 0$. Moreover, the multiplicity of each $E_n$ as an eigenvalue is $2n+1$ and the eigenvectors admit an explicit expression.

    For more general systems, the existence of such a ground state still holds in the case of positively charged or neutral atoms, that is, when $\sum_{m=1}^M Z_m \geq N$ (see \cite{hunziker2000quantum}), but explicit expressions for the ground energy and the ground state are no longer available. Nevertheless, these quantities are of particular interest for chemists, and a wide range of numerical methods have been developed in order to provide a numerical approximation of the ground state energy.

    Several families of methods exists, including wave-function methods, which are based on a direct approximation of the wave-function, Density Functional Theory~\cite{Kohn1965-xz,lieb2002density} based approximations which mainly aim at approximating the electronic density instead of the wavefunction, as well as Monte--Carlo methods which aims at minimizing the Rayleigh quotient $\frac{\langle \psi, H \psi \rangle}{\langle \psi, \psi \rangle}$ using Monte--Carlo sampling for the evaluation of the high-dimensional integrals.

    We will start here by presenting two wavefunction methods which consist in  choosing some particular subset of $\Sigma$ and minimize in~\eqref{equ: minmax principle} over this space instead of $\Sigma$. Note that it is crucial that the approximation satisfies the physical constraints given by the problem, that is antisymmetry.
    One way to ensure this is to base the approximation on so-called Slater determinants, defined for a given family of $N$ orthogonal functions $\varphi_1,\ldots, \varphi_N$ in $L^2(\mathbb{R}^3)$, as follows: for almost all $\bfm{r}_1, \cdots , \bfm{r}_N \in \mathbb{R}^3$,
 \begin{equation} \label{equ: Slater determinant}
                \varphi_1 \wedge \cdots \wedge \varphi_N (\bfm{r}_1, \cdots , \bfm{r}_N)
                :=   \frac{1}{\sqrt{N!}} \det(\varphi_j(\bfm{r}_i))_{i,j=1,...,N}.
            \end{equation}
    The approach which consists in choosing an orthonormal family of functions $(\varphi_k)_{k=1}^K$ of $L^2(\mathbb R^3)$ so that $\phi_k \in H^1(\mathbb{R}^3)$ for all $1\leq k \leq K$, and using the set
        \[
            \Sigma_K = \mathrm{Span} \, (\varphi_{j_1} \wedge ... \wedge \varphi_{j_N})_{1 \leq j_1 < j_2 < ... < j_N \leq K}
        \]
        leads to a discrete linear subspace $\Sigma_K$ the dimension of which rises exponentially with the number of electrons, making this method impossible to use even for seemingly small systems. This phenomenon is known as the curse of dimensionality, and is a general issue for the resolution of high dimensional problems.

\subsection{Hartree--Fock method}

In the well-known Hartree--Fock method~\cite{Lieb1977-ps}, an anti-symmetric wave-function is constructed by means of a single Slater determinant. Problem~\eqref{equ: minmax principle} is then replaced by
        \begin{equation} \label{equ: Hartee-Fock minmax}
            E_{HC}
            :=   \inf_{\psi \in \Sigma_{HF}, \norm \psi = 1} \langle \psi, H\psi \rangle,
        \end{equation}
        with
        \begin{equation}
            \Sigma_{\rm HF}
            =   {\rm Span} \left\{\varphi_1 \wedge ... \wedge \varphi_N, \quad 
            \varphi_1, \ldots, \varphi_N \in V, \quad \langle \varphi_i, \varphi_j\rangle = \delta_{ij} \text{ for } i,j=1,\ldots, N\right\},
        \end{equation}
        where $V$ is a finite dimensional subspace of $H^1(\R^3)$.
        Since $\Sigma_{\rm HF} \subset \Sigma$, it holds that $E_{\rm HC}\geq E_0$, but the equality does not hold in general.
        The space $\Sigma_{\rm HF}$ has a natural parametrization by collections of $N$ functions with variables in $\R^3$, hence the Hartree--Fock method does not suffer from the curse of dimensionality, and scales in general cubically with respect to the dimension of $V$.
        
        The approximation can be improved further by allowing two or more generally a linear combination of Slater determinants. This leads to the so-called Configuration Interaction (CI) methods, that we do not detail here further. A thorough presentation of classical methods in quantum chemistry can be found in \cite{Helgaker2014-py,cancesmadaylebris2003computational}.

\subsection{Density matrix renormalization group (DMRG)}

    One way to represent high dimensional vectors is to rely on tensor approximations, in particular the tensor train format (see \cite{dupuy2023lecture}). Sadly, these tensor representations are generally unable to deal directly with the antisymmetry constraint. One way to get rid of this constraint is to reformulate problem~\eqref{eq:schrod} using second quantization, which we present briefly.  We then give a short introduction to the tensor train format and the DMRG algorithm.
\subsubsection{Second quantization}
    Let us introduce the formalism of second quantization. For a more detailed description, see for example \cite{Helgaker2014-py,Legeza_2014}. The idea of second quantization is to fix some integer $d \in \mathbb{N}^*$ and family of functions $(\varphi_j)_{j=1}^d$ of $H^1(\R^3)$ and then define, for any sequence $\mu:=(\mu_j)_{j=1}^d \in \{0,1\}^{d}$, the Slater determinant:
    \begin{equation}
        \Psi_\mu
        :=   \bigwedge_{ j \text{ s.t.} \mu_j = 1 } \varphi_j,
    \end{equation}
    where $\bigwedge_{ j \text{ s.t.} \mu_j = 1 } \varphi_j$ corresponds to the Slater determinant defined in \eqref{equ: Slater determinant} taking the functions $\varphi_j$ such that the corresponding indices $\mu_j$ are equal to one. Notice that, for such a $\Psi_\mu$, $N_\mu:=\sum_{j=1}^\infty \mu_j$ is equal to the number of terms appearing in the Slater determinant $\Psi_\mu$, so that $\Psi_\mu\in \bigwedge_{i=1}^{N_\mu} L^2(\mathbb{R}^3)$.
    Taking an integer $d\in \mathbb{N}^*$, we then define the discrete Fock space of dimension $d$ as follows:
    \begin{equation} \label{equ: discrete Fock space}
        \mathcal F^d
        :=   \set{\sum_{\mu \in \{0,1\}^d} a_\mu \Psi_\mu}{(a_\mu)_{\mu \in \{0, 1\}^d} \in (\mathbb R^2)^{\otimes d}},
    \end{equation}
    where the sum must be understood in the sense of a formal linear combination. We observe in particular that if $N \leq d$, the following equality holds:
    \begin{equation} \label{equ: Fock space identification}
        \Sigma_N^d
        :=   \mathrm{Span} \, (\varphi_{j_1} \wedge ... \wedge \varphi_{j_N})_{1 \leq j_1 < j_2 < ... < j_N \leq d}
        =   \set{\sum_{\substack{\mu \in \{0,1\}^d \\ \sum_j \mu_j = N}} a_\mu \Psi_\mu}{(a_\mu)_{\mu \in \{0, 1\}^d} \in (\mathbb R^2)^{\otimes d}}
        \subset \mathcal F^d.
    \end{equation}
    
    Let $\mathcal T^d = (\mathbb R^2)^{\otimes d}$, and $(e_0, e_1)$ be the canonical basis of $\R^2$. Then a function $\mathcal J$ can be uniquely defined for any $\Psi_\mu$ (with $\mu=(\mu_1, \ldots, \mu_d) \in \{0, 1\}^d$) by
    \begin{equation} \label{equ: Fock tensor isomorphism}
        \mathcal J(\Psi_\mu)
        =   e_{\mu_1} \otimes e_{\mu_2} \otimes ... \otimes e_{\mu_d}.
    \end{equation}
    and extended by linearity to yield an isomorphism between $\mathcal T^d$ and $\mathcal F^d$.
    
    We also define the following matrices, belonging to $\mathbb{R}^{2\times 2}$,
    \begin{equation}
        A = \begin{pmatrix} 0 & 1 \\ 0 & 0\end{pmatrix}, \quad
        A^\dagger = \begin{pmatrix} 0 & 0 \\ 1 & 0\end{pmatrix}, \quad
        S = \begin{pmatrix} 1 & 0 \\ 0 & -1\end{pmatrix}, \quad
        I = \begin{pmatrix} 1 & 0 \\ 0 & 1\end{pmatrix},
    \end{equation}
    and, for $1 \leq p \leq d$, the operators acting upon $\mathcal T^d$:
    \begin{equation}
    \begin{array}{l}
        a_p = S \otimes ... \otimes S \otimes A \otimes I \otimes ... \otimes I
        \quad \text{where $A$ is in $p$-th position},\\
        a^\dagger_p = S \otimes ... \otimes S \otimes A^\dagger \otimes I \otimes ... \otimes I
        \quad \text{where $A^\dagger$ is in $p$-th position}.\\
        \end{array}
    \end{equation}
    
    The Hamiltonian $H$ can be represented as an operator acting on $\mathcal T^d$ as follows (see \cite[Theorem 2.1]{Legeza_2014} ).

    \begin{proposition}
        Define the operator $\tilde H = \mathcal J H \mathcal J^{-1}$ on $\mathcal T^d$. It holds that
            \begin{equation} \label{equ: Fock hamiltonian}
            \tilde H
            =   \sum_{p, q=1}^d h_{pq} a_p^\dagger a_q + \sum_{a, b, p, q=1}^d g_{abpq} a_a^\dagger a_b^\dagger a_p a_q,
        \end{equation}
        where
        \begin{equation*}
            h_{pq}
            =   \frac 1 2 \int_{\R^3} \overline{\nabla \varphi_q}(\bfm{r}) \cdot \nabla \varphi_p(\bfm{r}) d\bfm{r}
                - \int_{\R^3} \sum_{m=1}^M \frac{Z_m}{\abs{\bfm{r} - \bfm{X}_m}} \overline{\varphi_q}(\bfm{r}) \varphi_p(\bfm{r}) d\bfm{r},
        \end{equation*}
        and
        \begin{equation*}
            g_{abpq}
            =   \int_{\R^3 \times \R^3} \frac{\overline{\varphi_a}(\bfm{r}) \overline{\varphi_b}(\bfm{r}') \varphi_q(\bfm{r}) \varphi_p(\bfm{r}')}{\abs{\bfm{r} - \bfm{r}'}} d\bfm{r} d\bfm{r}'.
        \end{equation*}
    \end{proposition}
    
    Defining also
    \begin{equation}
        P = \sum_{p=1}^d a_p^\dagger a_p,
    \end{equation}
    it can be checked that $\Sigma_N^d$ is in fact equal to the space $\ker (\mathcal J P \mathcal J^{-1} - N)$.
    Indeed, for any $1 \leq p \leq d$, $a_p^\dagger a_p = I \otimes ... \otimes B \otimes I \otimes ... \otimes I$ where $B := \begin{pmatrix} 0 & 0 \\ 0 & 1 \end{pmatrix}$ is in $p$-th position, therefore for any $(a_\mu)_{\mu \in \{0, 1\}} \subset \mathbb R$,
    \[
    \begin{split}
        \mathcal J P \mathcal J^{-1} \left( \sum_{\mu \in \{0, 1\}^d} a_\mu \Psi_\mu \right)
        &=   \sum_{p=1}^d J a_p^\dagger a_p J^{-1} \sum_{\mu \in \{0, 1\}^d} a_\mu \Psi_\mu
        =  \sum_{\substack{p = 1, ..., d \\ \mu \in \{0, 1\} \\ \mu_p = 1}} a_\mu \Psi_\mu
        =   \sum_{\mu \in \{0, 1\}^d} \left( \sum_{p=1}^d \mu_p \right) a_\mu \Psi_\mu.
    \end{split}
    \]

    It follows that the minimization problem~\eqref{equ: minmax principle} restricted to the space $\Sigma_N^d$, i.e.
    \begin{equation}
        \inf_{\psi \in \Sigma_N^d, \, \norm{\psi}_{L^2} = 1} \langle \psi, H \psi \rangle
    \end{equation}
    is equivalent, if $N \leq d$, to
    \begin{equation} \label{equ: Fock minmax}
        \inf_{\substack{T \in \mathcal T^d, \, \norm T = 1 \\ PT = NT}} \langle T, \tilde H T \rangle,
    \end{equation}
    where the notation $\langle \cdot, \cdot \rangle$ and $\|\cdot\|$ refers here to the euclidean scalar product and norm of $\mathcal T^d$. 
    We therefore now need to solve a constrained optimization problem over the tensor space $\mathcal T^d = (\mathbb R^2)^{\otimes d}$ on which it is possible to use classical tensor formats, such as the tensor train format presented below, since the antisymmetry constraint has been removed.

\begin{remark}[Choice of basis functions]
In principle, one can choose any orthonormal basis of functions in $L^2(\mathbb R^3)$, such as finite elements or Hermite functions. However, using these functions does not lead in general to good results in practice. The reason is that the cost of the calculation increases at the fourth power with $d$. Indeed, the sum in~\eqref{equ: Fock hamiltonian} has $d^2+d^4$ terms.
Therefore conducting an efficient calculation requires to pick a set of orthogonal functions that is as small as possible. Chemists have dedicated a lot of effort in order to find such optimized basis sets, which are numerous nowadays and include for example  the Pople basis sets "STO-NG", or the correlation consistent basis sets "cc-pVDZ", "cc-pVTZ", etc. A monograph on the subject can be found in~\cite{Perlt2021-ya}. In Python for instance, these bases are implemented in the library PySCF~\cite{Sun2018-jt}. 
In the numerical tests with DMRG we restricted ourselves to a few relatively small basis sets. Indeed in order to keep the computational power used for the numerical experiments reasonable, we limited ourselves to about 10 functions per basis set.
\end{remark}

\subsubsection{Tensor Train (TT) format} \label{sec: Tensor Train format}

    The second quantization formalism reduces problem \eqref{equ: minmax principle} to another problem \eqref{equ: Fock minmax} on the tensor space $\mathcal T^d = (\mathbb R^2)^{\otimes d}$, hence removes the antisymmetry constraints, but does not reduce the dimensionality of the problem. Indeed, an element $(a_\mu) \in \mathcal T^d$ is still a tensor with $2^d$ coefficients and therefore often too large for its coefficients to be stored directly.
    The tensor train format offers a convenient way to approximate such elements of $\mathcal T^d$ without storing too many coefficients. We now present the main ideas of the tensor train format, a detailed presentation can be found in \cite{dupuy2023lecture}.

    The goal is to approximate the elements of the tensor space $\mathbb R^{n_1} \otimes ... \otimes \mathbb R^{n_d}$ for $n_1, n_2, ..., n_d \geq 1$, which can be seen as the space $ \mathbb R^{n_1 \times ... \times n_d}$, with elements $(S_{i_1, ..., i_d})_{1 \leq i_j \leq n_j}$.
    The idea of the tensor train format is to choose some positive integers $r_1, ..., r_{d-1}$, $r_0 = r_d = 1$, and matrices $C_j[i_j] \in \mathbb R^{r_{j-1} \times r_j}$ for $1 \leq j \leq d, 1 \leq i_j \leq n_j$ and define the associated tensor $S \in \R^{n_1 \times ... \times n_d}$ whose coefficients are given by
    \begin{equation} \label{equ: TT def formula}
        \forall (i_1, ..., i_d) \in \llbracket 1, n_1 \rrbracket \times ... \times \llbracket 1, n_d \rrbracket, \quad
        S_{i_1, ..., i_d} = \underbrace{C_1[i_1]}_{\in \mathbb R^{1 \times r_1}} \underbrace{C_2[i_2]}_{\in \mathbb R^{r_1 \times r_2}} \cdots \underbrace{C_d[i_d]}_{\in \mathbb R^{r_{d-1} \times 1}}.
    \end{equation}
    The matrices $C_j$ are called the cores of the tensor train, and the $r_j$ are called the ranks. The total cost to store the cores is $\sum_{j=1}^d n_j r_{j-1} r_j$ which is much smaller than the $n_1 n_2 ... n_d$ required to store the full tensor as long as the $r_j$ are not taken too large. Let us define $\tau$ as the $d$-linear application which maps the cores $(C_j)_{j=1}^d$ onto the tensor whose coefficients are defined by~\eqref{equ: TT def formula}.
    The advantage of the tensor train format is that the result of basic operations such as the inner product can be computed easily. For example, let $S = \tau(C_1, ..., C_d)$ and $T = \tau(D_1, ..., D_d)$), then
    \[
    \begin{split}
        \langle S, T \rangle
        &=  \sum_{i_1, ..., i_d} S_{i_1, ..., i_d} T_{i_1, ..., i_d}    \\
        &=  \sum_{i_1, ..., i_d} C_1[i_1] C_2[i_2] ... C_N[i_d] D_N[i_d]^T ... D_2[i_2]^T D_1[i_1]^T   \\
        &=  \sum_{i_1, ..., i_{d-1}} C_1[i_1] ... \left( \sum_{i_d} C_d[i_d] D_d[i_d]^T \right) ... D_1[i_1]^T.
    \end{split}
    \]
    Therefore, one can first compute the sum with respect to $i_N$, then the one with respect to $i_{N-1}$, and so on, which leads to much less operations than computing the sum over all multi-indices $(i_1, ...., i_d)$.

    A similar format, called the tensor train operator format (TTO), can be used to store linear operators over tensors :
    \begin{equation} \label{equ: TTO def formula}
        H_{i_1, ..., i_d}^{j_1, ..., j_d} = \underbrace{P_1[i_1, j_1]}_{\in \mathbb R^{1 \times r_1}} \underbrace{P_2[i_2, j_2]}_{\in \mathbb R^{r_1 \times r_2}} ... \underbrace{P_d[i_d, j_d]}_{\in \mathbb R^{r_{d-1} \times 1}}.
    \end{equation}
    Typically, an operator of the form
    \begin{equation} \label{equ: Example operator TTO}
        H = F_1 \otimes ... \otimes F_d + G_1 \otimes ... G_d,
    \end{equation}
    with $F_k \in \R^{n_k \times n_k}$ can be easily expressed as a TTO.
    To see this, we express \eqref{equ: Example operator TTO} coefficient-wise, which leads to
    \begin{equation*}
        H_{i1, ..., i_d}^{j_1, ..., j_d}
        =   (F_1)_{i_1}^{j_1} ... (F_d)_{i_d}^{j_d} + (G_1)_{i_1}^{j_1} ... (G_d)_{i_d}^{j_d}.
    \end{equation*}
    Therefore, \eqref{equ: TTO def formula} holds with
    \begin{equation*}
        P_1[i_1, j_1]
        =   \begin{pmatrix} (F_1)_{i_1}^{j_1} & (G_1)_{i_1}^{j_1} \end{pmatrix}
    \end{equation*}
    \begin{equation}
        P_k[i_k, j_k]
        \begin{pmatrix}
            (F_k)_{i_k}^{j_k} & 0 \\
            0 & (G_k)_{i_k}^{j_k}
        \end{pmatrix}
        , \quad
        \text{for $2 \leq k \leq d-1$},
    \end{equation}
    \begin{equation*}
        P_d[i_d, j_d]
        =   \begin{pmatrix} (F_d)_{i_d}^{j_d} \\ (G_d)_{i_d}^{j_d} \end{pmatrix}
    \end{equation*}    
    Similarly, a larger sum of tensor operators such as \eqref{equ: Fock hamiltonian} can also be expressed as a TTO.
    
    Another significant advantage of the TTO format is the compatibility with the TT format in the sense that we can apply a TTO to a TT and obtain the result as a TT, however with larger ranks. Explicitly if $S$ and $H$ are defined as in \eqref{equ: TT def formula} and \eqref{equ: TTO def formula}, and
    \[
        T_{i_1, ..., i_d}
        =   \sum_{j_1, ..., j_d} H_{i_1, ..., i_d}^{j_1, ..., j_d} S_{j_1, ..., j_d},
    \]
    then $T$ can be expressed as a TT with cores $D_1, ..., D_d$ defined by
    \[
        D_k[i_k]
        =   \sum_{j_k} P_k[i_k, j_k] \otimes C_k[j_k].
    \]

\subsubsection{DMRG algorithm}

    With all these operations defined on tensors, we are now able to write the DMRG algorithm, which will be used to compute an approximation of the ground state of \eqref{equ: Fock hamiltonian}.
    We first state the algorithm in a general abstract setting. Let $J$ be a functional
    \[
        J : \mathbb R^{n_1 \times ... \times n_d} \to \mathbb R
    \]
    that we want to minimize. Then we can apply the following procedure:
    \begin{algorithm}[H]
        \caption{DMRG (simplified)} \label{algo: DMRG}
        \begin{algorithmic}[1]
            \REQUIRE Pick some initial tensor cores $C_1^0, ..., C_d^0$ randomly
            \REQUIRE $n = 0$
            \WHILE{not converged}
                \FOR{$i = 1$ to $d-1$}
                    \STATE $C_i^{n+1} \leftarrow \underset{C_i}{\text{argmin}} \, J \circ \tau(C_1^{n+1}, ..., C_{i-1}^{n+1}, C_i, C_{i+1}^n, ..., C_d^n)$
                \ENDFOR
                \FOR{$i = d$ to $2$}
                    \STATE $C_i^{n+1} \leftarrow \underset{C_i}{\text{argmin}} \, J \circ \tau(C_1^{n+1}, ..., C_{i-1}^{n+1}, C_i, C_{i+1}^n, ..., C_d^n)$
                \ENDFOR
                \STATE $n \leftarrow n + 1$
            \ENDWHILE
        \end{algorithmic}
    \end{algorithm}
    In practice, the DMRG algorithm also involves some orthogonality constraints on the cores, but we do not give the details here for simplicity (see \cite{dupuy2023lecture}).
    The most important observation is that at each iteration, one only solves an optimization problem with respect to a single core which corresponds to a simple quadratic problem and a limited number of parameters.

    Going back to the Fock Hamiltonian~\eqref{equ: Fock hamiltonian} and problem \eqref{equ: Fock minmax}, the functional $J$ is taken in our case  as
    \[
        J(T)
        =   \frac{\langle T, \tilde H T \rangle}{\langle T, T \rangle}.
    \]
    As we saw in Section \ref{sec: Tensor Train format}, the operator $\tilde H$ can easily be expressed as a TTO, and the inner product of two TTs can be computed. Therefore, if $T$ is given as a TT, computing $J(T)$ only requires a number of operations which scales linearly with $d$. For this reason, Algorithm \ref{algo: DMRG} does not suffer from the curse of dimensionality.

\subsection{Neural networks based algorithms}
Neural networks have emerged as powerful tools to approximate wavefunctions and have provided a way to overcome the computational expensiveness of highly accurate methods with respect to the number of electrons. 
We present here two of the most emerging methods (PauliNet \cite{Hermann_2020} and FermiNet \cite{Pfau_2020}) that have combined quantum Monte-Carlo and deep learning strategies, and have been recognized as providing a good trade-off between accuracy and computational cost.  Unlike DMRG, these methods work directly with the wavefunctions and do not require the formalism of second quantization. Our goal is to describe them in the simplest possible way and as precisely as possible. We start by a presentation of the variational Monte--Carlo method, before presenting the functional form used in the neural network based on the Slater Jastrow-backflow ansatz. Finally, we present FermiNet and PauliNet in a unified way.

\subsubsection{Variational Monte-Carlo}

The deep neural networks (DNNs) that we study here are trained using variational quantum Monte--Carlo (VMC) methods. 
The minimization problem~\eqref{equ: minmax principle} that is approximately solved using a neural network may exhibit convergence issues due to the  constraint on the $L^2$ norm of the wavefunction. Therefore, the constraint is replaced by the evaluation of the Rayleigh quotient of the form, for a given wavefunction $\psi_{\theta}$ depending on parameters $\theta$
\begin{align*}
    E[\psi_{\theta}] = \frac{\langle \psi_{\theta}, H \psi_{\theta}\rangle }{\langle \psi_{\theta},\psi_{\theta} \rangle}.
\end{align*}
Moreover, this energy can be rewritten as follows
\[
E[\psi_{\theta}] = \mathbb{E}_{x \sim p}[\psi_\theta^{-1}(\mathbf{r})H\psi_\theta( \bfm{r})]
\]
where $p(  \bfm{r}) \sim \psi_\theta^2(  \bfm{r})$ is a probability measure.
The wavefunction is then represented by a neural network with weights $\theta$ that are optimized to minimize the energy. 
In practice, the computation of the energy 
involves the computation of high-dimensional integral that are approximated by sampling from probability distribution $ p(\bfm{r}) $ using Markov Chain Monte Carlo (MCMC) techniques, and the energy is then approximated by a classic Monte--Carlo sum.
\subsubsection{Neural networks}
Neural networks (NNs) incredible capacity of approximating functions are no longer to prove. The algorithms we study here mostly use feed-forward neural networks, which can be represented as nonlinear parametric functions $ u_{\theta} : \mathbb R^{d_{in}} \rightarrow \mathbb R^{d_{out}}$ where $u_{\theta}$ can be written as : 
\begin{align*}
    u_{\theta}(x_0) := \sigma_L(W_L x_{L-1} + b_L) \circ \sigma_{L-1}(W_{L-1} x_{L-2} + b_{L-1}) \circ \ldots \circ \sigma_1(W_1 x_0 + b_1)
\end{align*}
$ x_0 \in \mathbb R^{d_{in}}$ represents the input vector of the neural network, the $ (\sigma_\ell)_{1 \leq \ell \leq L}$ are continuous non-polynomial functions called \textit{activation} functions, the $ \sigma(W_\ell x_{\ell - 1} + b_\ell) \enspace, 1 \leq \ell \leq L - 1 $ are the \textit{hidden layers}, and $ \sigma_L(W_L x_{L - 1} + b_L)$ is the output layer, with $ L $ being the number of layers of the NN. The weights $ (W_\ell)_{1 \leq \ell \leq L} $ and bias $ (b_\ell)_{1 \leq \ell \leq L} $ are the parameters to be optimized. We note $ \theta \in \mathbb R^p $ the set of parameters.

Training a neural network typically involves solving  a minimization problem
\begin{align*}
    \text{Find } \theta^* \in \mathbb R^p \enspace | \enspace \theta^* := \underset{\theta \in \mathbb R^p}{\text{argmin }} C(\theta)
\end{align*}
where $ C $ is a cost function. 

In the context of the electronic Schrödinger equation, the cost function can be set to be the energy we want to minimize, and therefore the problem becomes
\begin{align*}
    \text{Find } \theta^* \in \mathbb R^p \enspace | \enspace \theta^* := \underset{\theta \in \mathbb R^p}{\text{argmin }} E[\psi_{\theta}]
\end{align*}
This is an unsupervised learning problem, as the cost function is only dependent on the input data and the network's output (the wavefunction).
\subsubsection{The Slater Jastrow-Backflow ansatz}

The Slater Jastrow-Backflow ansatz is a powerful method that aims at better capturing the correlation effects present in the system, to accurately represent the wavefunction. It is composed of three main ingredients. 

First, a basic component is the Slater determinant, described in~\eqref{equ: Slater determinant}. This ensures that the wavefunction satisfies the Pauli principle. The second component is the Jastrow factor, which consists in multiplying the Slater determinant by $e^{J(\bfm{r}_1, \ldots, \bfm{r}_N) }$, where $ J(\bfm{r}_1, \ldots, \bfm{r}_N) = - \sum_{i<j} U(\lvert \bfm{r}_i - \bfm{r}_j \rvert )$ is the Jastrow factor, with $ U $ a correlation function usually chosen as an exponential or polynomial function that decays with increasing particle separation, reflecting the tendency of electrons to repel each other at short distances. Note that this Jastrow factor is symmetric with respect to the permutation of variables and hence does not break the antisymmetry of the Slater determinant when multiplied by a Jastrow factor.

The third component is the backflow transform, which makes the orbitals in the Slater determinant dependent on the coordinates of all the electrons. This helps incorporating additional correlation effects beyond those captured by the Slater determinant and the Jastrow factor, therefore further improving the accuracy of wavefunction description. The backflow transformation can be represented as
    \begin{align}
        \phi_i(\bfm{r}_j) \rightarrow \phi_i(\bar{\bfm{r}_j}) \text{ where } \bar{\bfm{r}_j} = \bfm{r}_j + \sum_{k \not= j} \eta(\lvert \bfm{r}_j - \bfm{r}_k \rvert) \bfm{r}_k
        \label{backflow}
    \end{align}
    where $ \eta $ is a correlation function. 

We now present the two neural network based methods FermiNet and PauliNet, which take ideas from the Slater Jastrow backflow ansatz, and further improve it by using neural networks. 
\subsubsection{Fermionic Neural Networks}

    Fermionic neural networks (FermiNet~\cite{Spencer2021-ee})  is inspired from the Slater Jastrow backflow ansatz. The key idea behind the model is to allow the orbitals to depend on all electron positions, without explicitly writing the dependence. This allows a more flexible architecture, as the network can in principle learn hidden correlations between electrons. 

Denoting by $ \bfm{r}_{\not=j} = (\bfm{r}_1, \ldots, \bfm{r}_{j-1}, \bfm{r}_{j+1}, \ldots, \bfm{r}_N)$, the orbitals are defined as  functions depending on all variables, that is $ \varphi_i(\bfm{r}_j, \bfm{r}_{\not=j}) $, with the requirement that they are invariant with respect to the permutations of $ (\bfm{r}_{\not=j})$, i.e.
\begin{align}\label{eq:equivariance}
   \forall s \in S_N^{\not=j}, \, \, \, \,  \varphi_i(\bfm{r}_j, \bfm{r}_{\not=j}) = \varphi_i(\bfm{r}_j, \bfm{r}_{s(1)}, \ldots, \bfm{r}_{s(j-1)}, \bfm{r}_{s(j+1)}, \ldots, \bfm{r}_{s(N)})
\end{align}
where $ S_N^{\not=j} $ is the set of permutations of $\{ 1, \cdots, j-1, j+1, \cdots ,N\}$. This ensures that a determinant defined by the orbitals $\varphi_1,\ldots, \varphi_N$ as
\[
\begin{vmatrix}
     \varphi_1(\bfm{r}_1, \bfm{r}_{\not=1}) & \ldots & \varphi_1(\bfm{r}_N, \bfm{r}_{\not=N}) \\ 
    \vdots & \vdots & \vdots \\ 
    \varphi_N(\bfm{r}_1, \bfm{r}_{\not=1}) & \ldots & \varphi_N(\bfm{r}_N, \bfm{r}_{\not=N})
    \end{vmatrix}
\]
satisfies the antisymmetry constraint.
    
The wavefunction is then expressed as a linear combination of such determinants, namely
\begin{align}
    \psi(\bfm{r}_1, \ldots, \bfm{r}_N) = 
    \sum_{p=1}^P \lambda_p
    \begin{vmatrix}
     \varphi_1^k(\bfm{r}_1, \bfm{r}_{\not=1}) & \ldots & \varphi_1^k(\bfm{r}_N, \bfm{r}_{\not=N}) \\ 
    \vdots & \vdots & \vdots \\ 
    \varphi_N^k(\bfm{r}_1, \bfm{r}_{\not=1}) & \ldots & \varphi_N^k(\bfm{r}_N, \bfm{r}_{\not=N})
    \end{vmatrix},
\end{align}
where the $ (\lambda_k)_k $ are fixed weights.

The complete architecture of FermiNet is given in Algorithm \ref{algo:ferminet} where $L$ denotes the number of layers in the neural network, $P$ is the number of determinants in the linear combination, $\sigma $ is the nonlinear activation function. We denote in red the trainable parameters.

    \begin{algorithm}
    \caption{FermiNet scheme}
    \label{algo:ferminet}
    \begin{algorithmic}[1]

    \STATE \textbf{Input}: Electron positions $ \bfm{\bfm{r}} = (\bfm{r}_1, \ldots, \bfm{r}_N )$, nuclear positions $\{\bfm X_\ell \}_{\ell \leq 1}^M $

    \STATE \textbf{Output}: Wavefunction $ \psi(\bfm{r}) $

    \FOR{electron $i = 1$ to $N$}
        \STATE $ X_i^0 \leftarrow \text{concatenate}(\bfm{r}_i - \bfm{X}_k, \lvert \bfm{r}_i - \bfm{X}_k \rvert, \forall 1 \leq k \leq M) \rightsquigarrow$  Electron-Nuclei interactions
        \STATE $ X_{ij}^0 \leftarrow \text{concatenate}(\bfm{r}_i - \bfm{r}_j, \lvert \bfm{r}_i - \bfm{r}_j \rvert, \forall 1 \leq j\leq N) \rightsquigarrow $ Electron-Electron interactions
    \ENDFOR
    \FOR{layer $\ell = 0$ to $L-1$}
    \STATE $ Y^\ell \leftarrow \frac{1}{N} \sum_{i=1}^N X_i^\ell \rightsquigarrow $ ensure invariance w.r.t electrons permutation
    \FOR{electron $i = 1$ to $N$}
    \STATE $ Y_i^\ell \leftarrow \frac{1}{N} \sum_{i=1}^N X_{ij}^\ell \rightsquigarrow$ ensure invariance w.r.t electrons permutation
    \STATE $ Z_i^\ell \leftarrow \text{concatenate}(X_i^\ell, Y^\ell, Y_i^\ell) \rightsquigarrow$ feature vector for electron $ i $ at layer $ \ell$ 
    \STATE $ X_i^{\ell + 1} \leftarrow \sigma(\red{W_1^\ell} Z_i^\ell  + \red{b_1^\ell}) + X_i^{\ell}$
    \STATE $ X_{ij}^{\ell + 1} \leftarrow \sigma(\red{W_2^\ell} Z_i^\ell  + \red{b_2^\ell}) + X_{ij}^{\ell}$
    
    \ENDFOR
    \ENDFOR
    \FOR{determinant $ p=1 $ to $P$}
    \FOR{orbital $ i=1 $ to $N$}
    \FOR{electron $ j=1 $ to $N$}
    \STATE $ e \leftarrow \sum_{k=1}^M \red{w_{ik}^p} \text{exp}(-\lvert \red{A_{ik}^p} (\bfm{r}_j - \bfm{X}_k)\rvert) \rightsquigarrow $ encodes the behavior of wavefunction when electrons are far away from nuclei
    \STATE $ \phi_i^p(\bfm{r}_j, \bfm{r}_{\not=j}) = (\red{W_i^p} \cdot X_j^L + \red{b_i^p})e \rightsquigarrow $ final layer to calculate the orbitals
    \ENDFOR
    \ENDFOR
    \ENDFOR
    \STATE $ \psi(r) \leftarrow \sum_{p} \lambda_p det[(\phi_i^p(\bfm{r}_j, \bfm{r}_{\not=j}))_{i, j}]$
    \end{algorithmic}
\end{algorithm}

To summarize, FermiNet takes as input electron-electron and electron-nucleus features. Besides the differences $ \bfm{r}_i - \bfm{X}_k $ and $ \bfm{r}_i - \bfm{r}_j $, the distances $ \lvert \bfm{r}_i - \bfm{X}_k \rvert $ and $ \lvert \bfm{r}_i - \bfm{r}_j \rvert $ are also added. This is helpful for learning a Jastrow factor and backflow, since the algorithm does not assume a closed form of those functions. These features are then fed to several layers (linear transformations with activation function) with residual connections. 
These layers respect by construction the constraint of invariance with respect to electrons permutation~\eqref{eq:equivariance}, as they take as input the sum over electrons. A final linear transformation is then applied to the output of the last intermediate layer, to approximate the mutli-electron orbitals for each determinant. These orbitals represent the output of the overall Fermionic Neural Network, and are used to evaluuate the wavefunction as the weighted sum of determinants.

\subsubsection{PauliNet}

PauliNet~\cite{Hermann_2020} is another deep learning approach using neural networks as the wavefunction ansatz. It also takes ideas from the Slater Jastrow backflow ansatz, but differs from FermiNet in the fact that it specifically encodes the Jastrow factor and backflow transform as trainable deep neural networks. More precisely, we define $ J_{\theta_1}$ and $ f_{\theta_2}$, respectively the neural networks parameterizing the Jastrow factor and backflow transform. Then, the wavefunction writes

    \begin{align}
        \psi_{\theta}(\bfm{r}) = e^{[c(\bfm{r}) + J_{\theta_1}(\bfm{r})]} \sum_{p=1}^P det[(\tilde{\varphi}_{\theta_2, p}^{j, i} (\bfm{r}))]
        \label{wavefunctionpn}
    \end{align}
    with 
    \begin{itemize}
        \item $ \tilde{\varphi}_{\theta_2, p}^{j, i}(\bfm{r}) = \varphi^j(\bfm{r}_i) f_{\theta_2, i}(\bfm{r}), \enspace \bfm{r} \in \mathbb R^{3N}, \bfm{r_i} \in \mathbb R^3 $, where $ (\varphi^j(\bfm{r}_i))_{i, j} $ are the one-electron molecular orbitals, $f_{\theta_2, i}(\bfm{r})$ is a backflow vector for electron $ i $ ($f_{\theta_2, i}(\bfm{r}) \in \mathbb R^N $ ), $ i $ represents the indices for electrons, $ j $ for the orbitals, and $p$ for the determinants.
        \quad
        \item $ c(\bfm{r}) = \sum_{i < j} \frac{c_{ij}}{1 + \lvert \bfm{r}_i - \bfm{r}_j \rvert} $ where $c_{ij}$ are coefficients related to spins.
    \end{itemize}

    For the wavefunction to be antisymmetric in this case, we need certain conditions on the neural networks. 
    \begin{itemize}
        \item $ J_{\theta} $ needs to be invariant with respect to permutation $ s_{ij} $ of electrons $i$ and $j$, i.e 
        \begin{align*}
            J_{\theta}(s_{ij}(\bfm{r})) = J_{\theta}(\bfm{r}), 
        \end{align*}
        \item $ f_{\theta}$ should be equivariant with respect to permutation $ s_{ij} $ of electrons $i$ and $j$, i.e.
        \begin{align*}
             s_{ij}f_{\theta, i}(\bfm{r}) = f_{\theta, j}(s_{ij}(\bfm{r})).
        \end{align*}
It can be easily checked that the overall antisymmetry of the wavefunction is satisfied when these two conditions are met.
    \end{itemize}

    We detail in Algorithm \ref{algo:paulinet} the explicit construction of the wavefunction.

    \begin{algorithm}[H]
    \caption{PauliNet scheme}
    \begin{algorithmic}[1]
    \STATE \textbf{Input} : Electron positions $ \bfm{r} $, nuclear positions $ \bfm{X} $
    \STATE \textbf{Output} : Wavefunction $ \psi(\bfm{r}) $
    \FOR{$i = 1$ to $N$}
        \STATE Initialization : $ x_i^{(0)} \in \mathbb{R}^{N \times d_e}$ (random trainable vector), 
        \FOR{$j=0$ to $ L - 1 $}
        \STATE $ x_i^{(j+1)} \leftarrow x_i^{(j)} + \red{F_{\theta_3}^{(j)}}(\{ x_k^{(j)}, \{ \lvert \bfm{r}_k - \bfm{r}_{k^\prime} \rvert\} \})$, where $ F_{\theta_3}$ is a graph convolutional neural network, equivariant w.r.t particle exchange.
    \ENDFOR
            \ENDFOR
    \STATE Compute single-electron orbitals $ \varphi^j (\bfm{r}_i) \enspace 1 \leq i, j \leq N $
    \STATE Define neural network $ \red{J_{\theta_1}} : \mathbb R^{d_e} \rightarrow \mathbb R $ with input $ \{ \sum_{i=1}^N x_i^L \}$
    \STATE Define neural network $ \red{f_{\theta_2, i}} : \mathbb R^{N \times d_e} \rightarrow \mathbb R^N $ with input $ \{ x_i^L\} $
    \STATE Define $ c(\bfm{r}) = \sum_{i < j} \frac{c_{ij}}{1 + \lvert \bfm{r}_i - \bfm{r}_j \rvert}$
    \STATE Calculate $ \psi_{\theta}(\bfm{r}) $ using formula (\ref{wavefunctionpn})
    \end{algorithmic}
    \label{algo:paulinet}    
\end{algorithm}
To summarize, PauliNet takes as input electrons and nuclear coordinates. These coordinates are first fed to a graph convolutional deep neural network to encode complex interactions between electrons. The output of this graph is a representation of each particle in a high dimension space ($ \mathbb R^{d_e}$). This output is then used as input for two neural networks, one parameterizing the Jastrow factor, and the other parameterizing the backflow transform. 
To guarantee that the Jastrow factor is invariant with respect to the permutation of the electrons, a sum over all electrons is taken as input to the neural network.
The backflow vectors $ f_{\theta_2, i} $ are neural networks with shared parameters, which take an input $ x_i^L $, which is equivariant with respect to electron permutation by construction of the graph. This therefore makes $ f_{\theta_2, i}$ equivariant as well. In addition, a multireference Hartree--Fock method is used to compute the single electron orbitals for a chosen number of determinants. Those orbitals helps starting the optimization in a satisfying region of convergence, and are not changed during the training process.

The graph neural network used in PauliNet is a modified version of the SchNet~\cite{schnet}, which in contrary to the initial algorithm where only nuclear coordinates and charges are used to model the energy, also takes into account the electron coordinates. 

\section{Numerical results} 
\label{sec:num}

We now present numerical results on the DMRG method as well as neural networks FermiNet and PauliNet methods. We start by providing a few convergence tests independently on these methods before comparing them on a few atoms and molecules.

\subsection{DMRG numerical experiments}

We applied the DMRG algorithm to the electronic eigenvalue problem formulated in the Fock space \eqref{equ: Fock minmax} on several atoms and small molecules (see Table~\ref{tab:tab1}). 
Although the DMRG algorithm breaks the curse of dimensionality in principle, the number of terms involved in the calculation of the second quantization Hamiltonian \eqref{equ: Fock hamiltonian} tends to increase rather drastically with the number of basis function, therefore making the DMRG computation expensive. Due to this limitation and the limited computational power available to us during CEMRACS where these numerical simulations were carried out, the number of electrons involved in our tests for DMRG could hardly exceed 10. 

\begin{figure}
        \centering
        \includegraphics[width=0.6\textwidth]{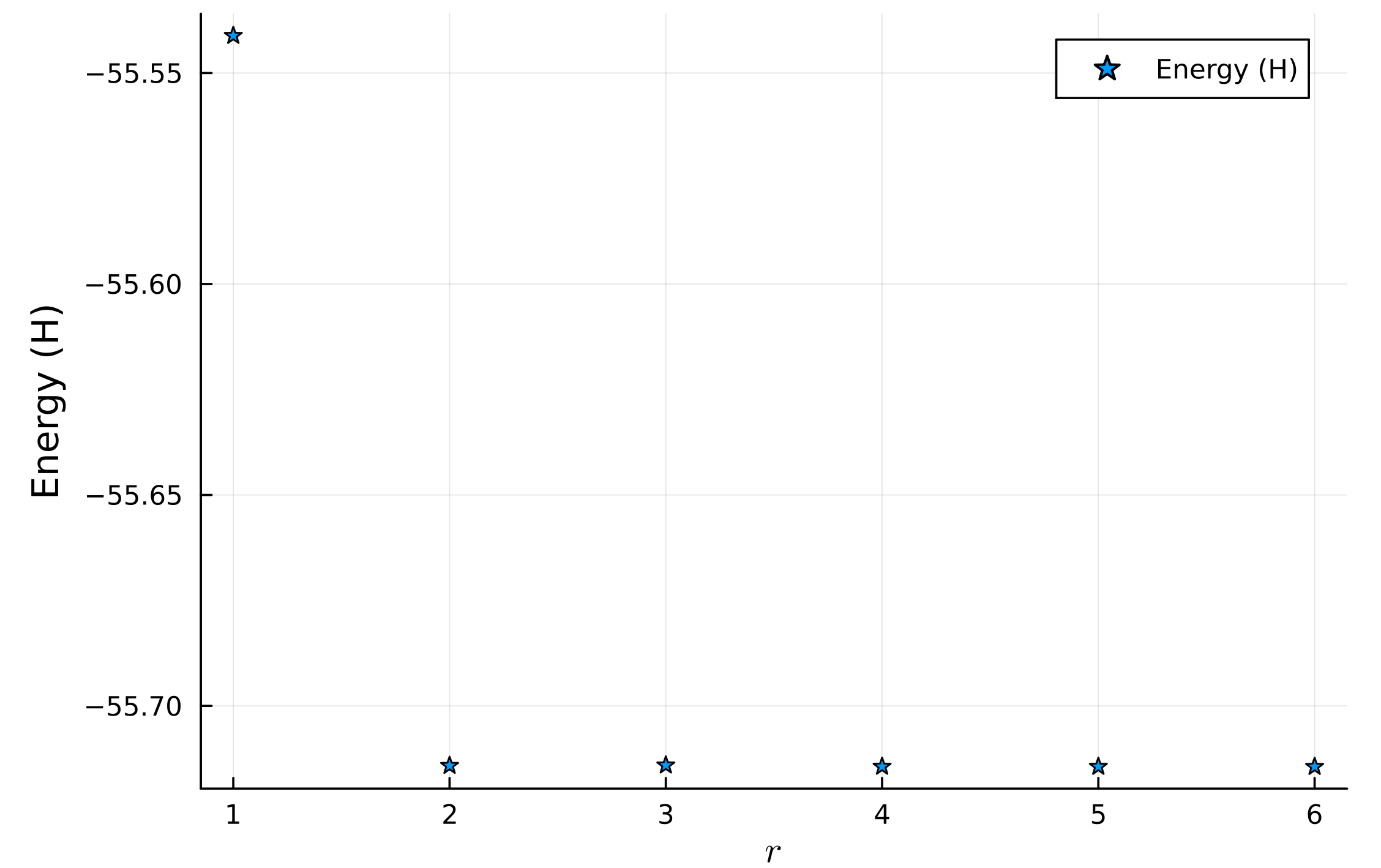}
        \caption{Energy (expressed in Hartree) computed by DMRG with the ranks $r_j$ of the tensor train all equal to $r$, with the basis "6-31G" for the molecule $LiH$}
        \label{fig: rank comparison}
    \end{figure}

\begin{figure}
        \centering
        \includegraphics[width=0.6\textwidth]{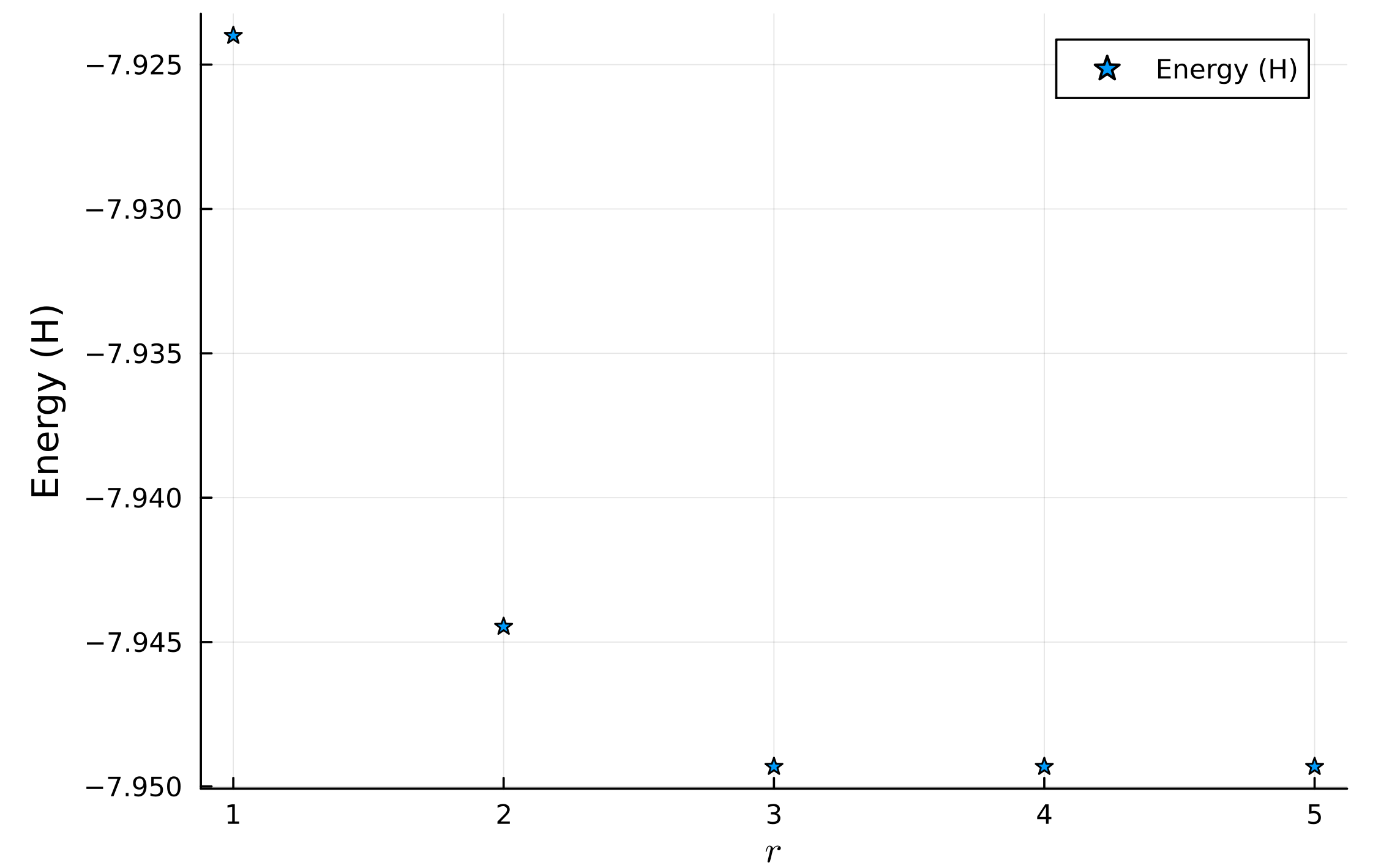}
        \caption{Energy (expressed in Hartree) computed by DMRG with the ranks $r_j$ of the tensor train all equal to $r$, with the basis "STO-6G" for the molecule $NH_3$}
        \label{fig: rank comparison 2}
    \end{figure}

We first tested the influence of the maximal chosen ranks defining the TT. We chose the ranks $r_1, ..., r_d$ of the tensor train such that they were all equal to $r$. In all our tests, we observed that the energy computed by DMRG quickly stabilises when the rank $r$ of the tensor train reaches a certain threshold, usually smaller than $6$. We illustrate this phenomenon on the $LiH$ molecule simulated with the chemical basis "6-31G",, which is a standard Pople basis set~\cite{Ditchfield1971-oq}, in Figure \ref{fig: rank comparison}, and the $NH_3$ molecule simulated with the chemical basis "sto-6g" in Figure \ref{fig: rank comparison 2}

\subsection{Numerical experiments on DNNs methods}

    We now present some numerical results related to the study of FermiNet and PauliNet. For all tests below, the optimizers used are the Kronecker Factorized Approximate Curvature (KFAC) \cite{kfac} for FermiNet, and Adam optimizer for PauliNet \cite{adam}. We use ethylene $ (C_2H_4) $ and ammonia $(NH_3)$ for the different experiments, and later provide the energy estimation for different atoms and molecules.

\subsubsection{Randomness}
On the one hand, aside from antisymmetry, there are no physical information encoded in FermiNet. The neural network is responsible for learning the orbitals, with adequate Jastrow factor and backflow transform. One the other hand, the network PauliNet specifically encodes the Jastrow and backflow as deep neural networks, and include more physics via electronic cusps and single-electron orbitals. 
For this reason, the ansatz in PauliNet gives less degrees of freedom than for FermiNet. 

\begin{figure}[H]
\minipage{0.5\textwidth}
    \includegraphics[width=\textwidth]{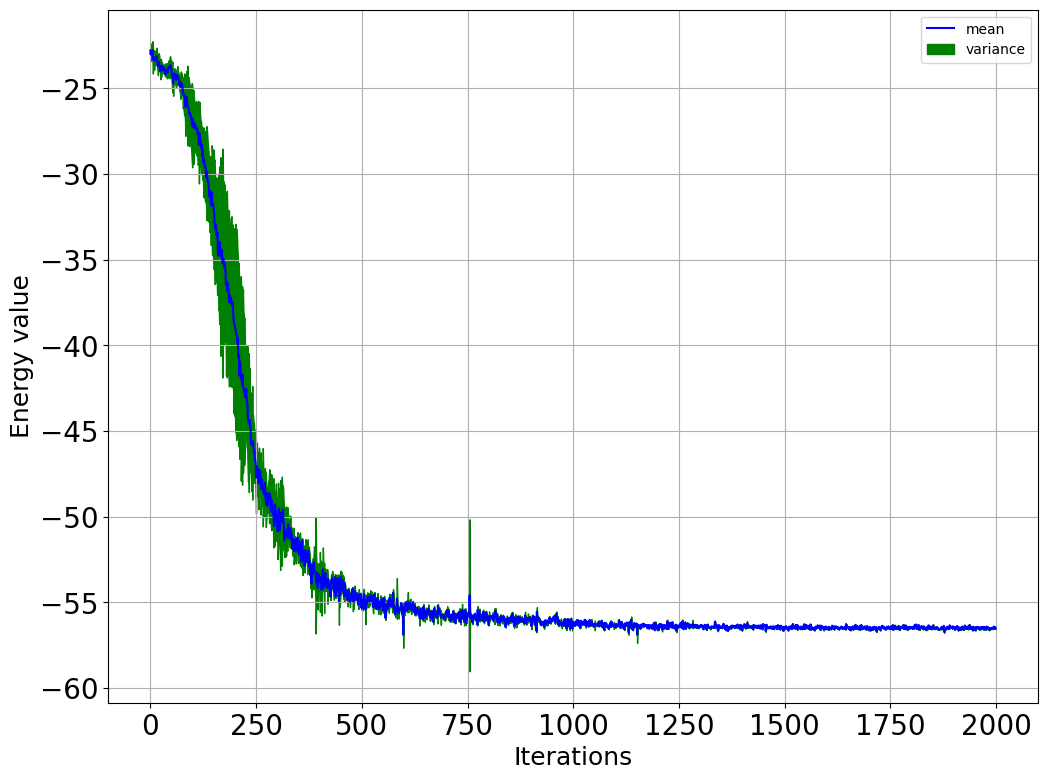}
        \caption{ Mean and variance for energy minimization of $NH_3$ with FermiNet}
        \label{fig:mean_variance_fn_nh3}
\endminipage
\minipage{0.5\textwidth}
\includegraphics[width=\textwidth]{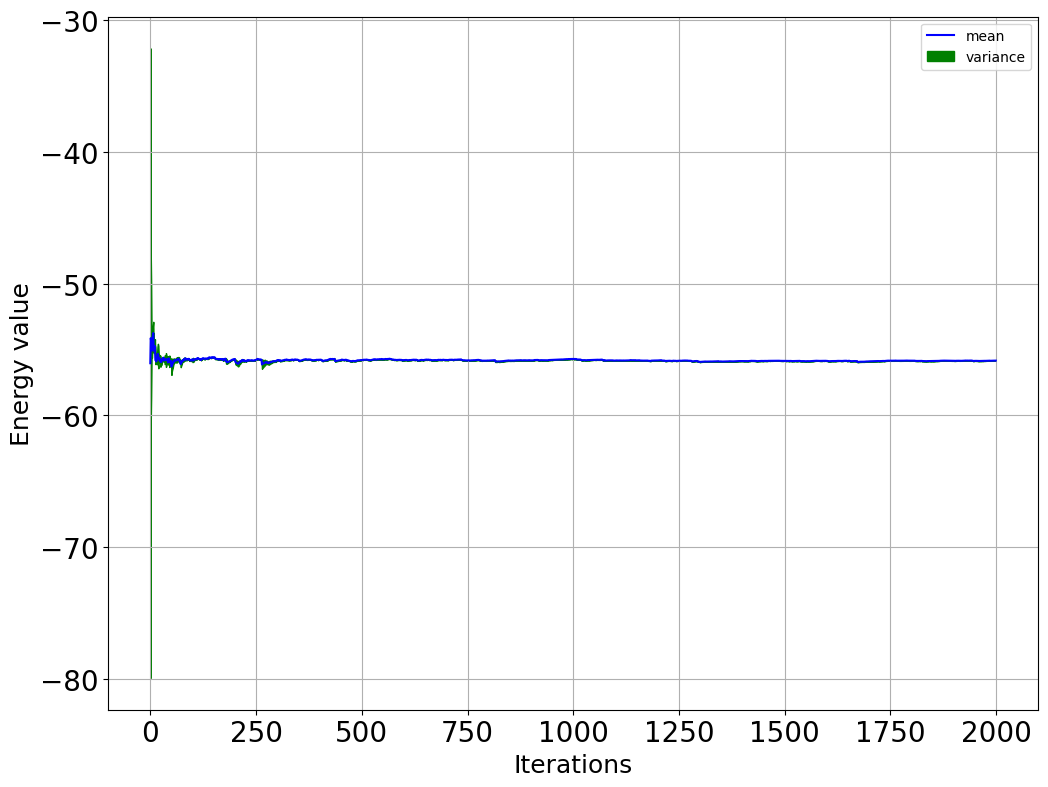}
    \caption{ Mean and variance for energy minimization of $NH_3$ with PauliNet}
    \label{fig:mean_variance_pn_nh3}
\endminipage
\end{figure}

\begin{figure}[H]
\minipage{0.5\textwidth}
    \includegraphics[width=\textwidth]{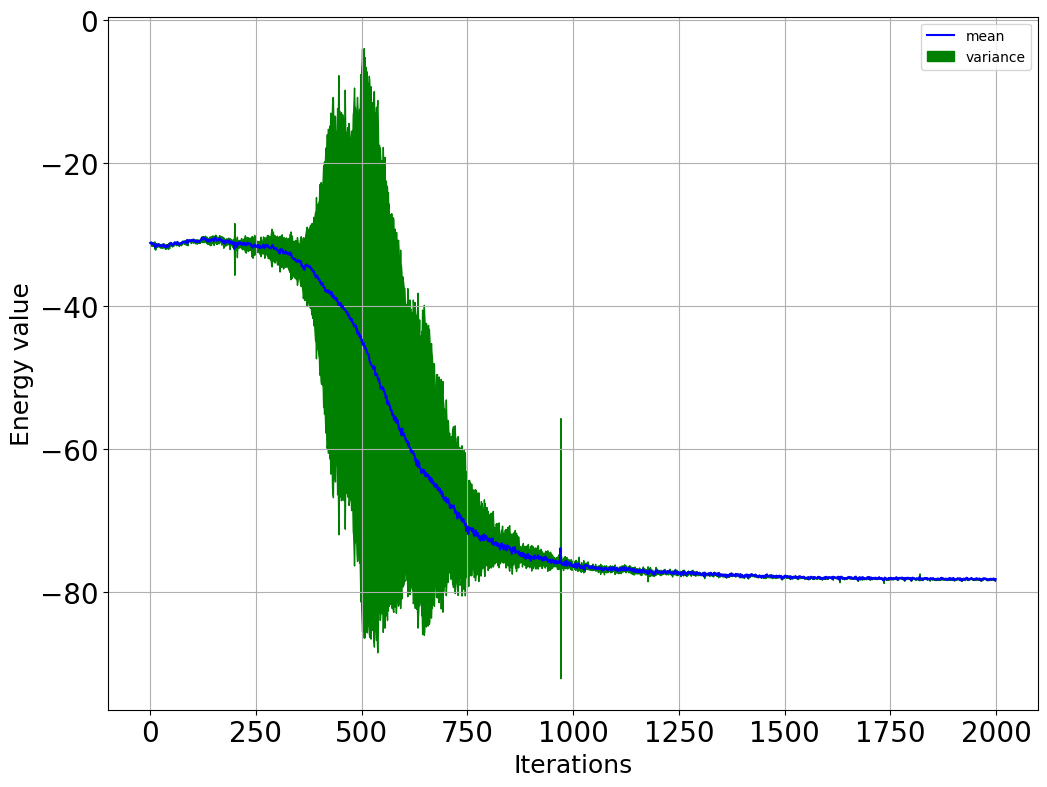}
        \caption{ Mean and variance for energy minimization of $C_2H_4$ with FermiNet}
        \label{fig:mean_variance_fn_c2h4}
\endminipage
\minipage{0.5\textwidth}
\includegraphics[width=\textwidth]{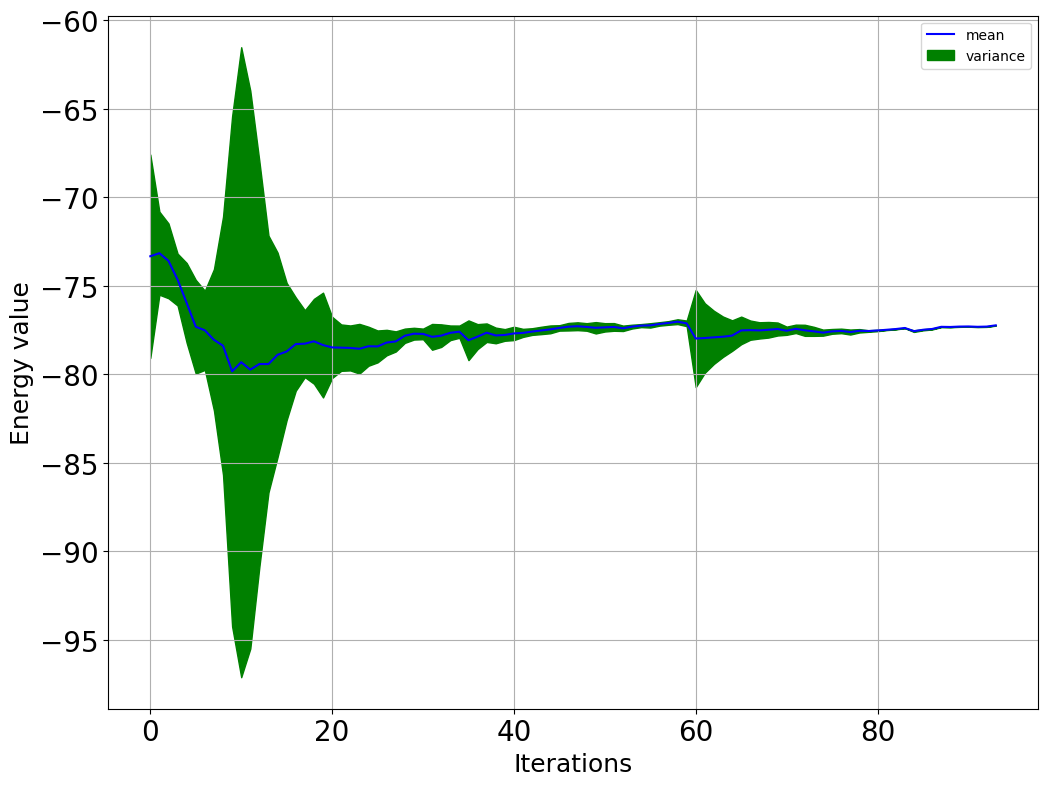}
    \caption{ Mean and variance for energy minimization of $C_2H_4$ with PauliNet}
    \label{fig:mean_variance_pn_c2h4}
\endminipage
\end{figure}

A natural question arising when dealing with neural networks is the stability of the model, how the model performs and behaves over the iterations of the optimization algorithm under different conditions. We investigate here, how different initializations impact the learning process in PauliNet and FermiNet, by performing multiple trainings. More specifically, we test the robustness of the models here $ NH_3 $ and $ C_2H_4$, by training them $ 10 $ times. We then display the mean curve and variance between the different learning processes. The results (Figure (\ref{fig:mean_variance_fn_nh3}),  Figure (\ref{fig:mean_variance_pn_nh3}), Figure (\ref{fig:mean_variance_fn_c2h4}), Figure (\ref{fig:mean_variance_pn_c2h4})) show tremendous stability of PauliNet for $ NH_3 $ as there is almost no variance observed during the training phase. Most importantly, both models appear to be robust, as the learned ground state is always the same.

\subsubsection{Pretraining step}

FermiNet is usually pretrained using Hartree--Fock orbitals, in order to achieve faster convergence. Roughly speaking, this means that we first train the network to solve a regression problem where the target is the eigenstate given by Hartree--Fock, allowing the weights during the training process, to be initialized in a region of faster convergence.
Here we train the model to see how much the pretraining step helps to converge faster towards the ground state. The model was pretrained using the basis set "STO-3G" and Adam optimizer. 

\begin{figure}[H]
\minipage{0.5\textwidth}
    \includegraphics[width=\textwidth]{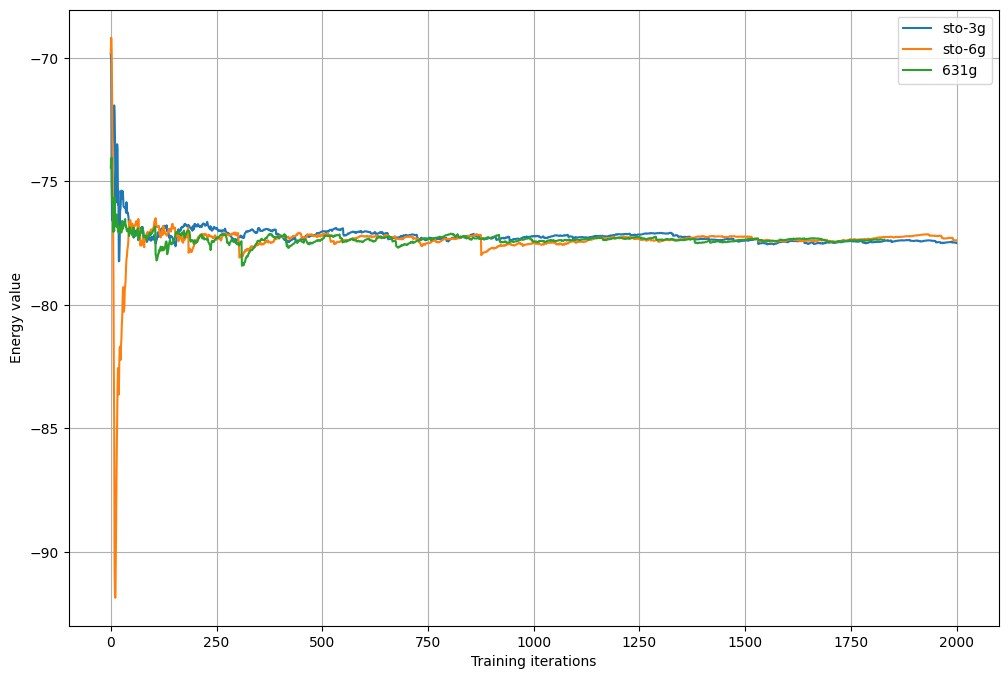}
        \caption{Paulinet : Energy minimization of $C_2H_4$ with different basis sets}
        \label{fig:pretraining_pn_c2h4}
\endminipage
\minipage{0.5\textwidth}
\includegraphics[width=\textwidth]{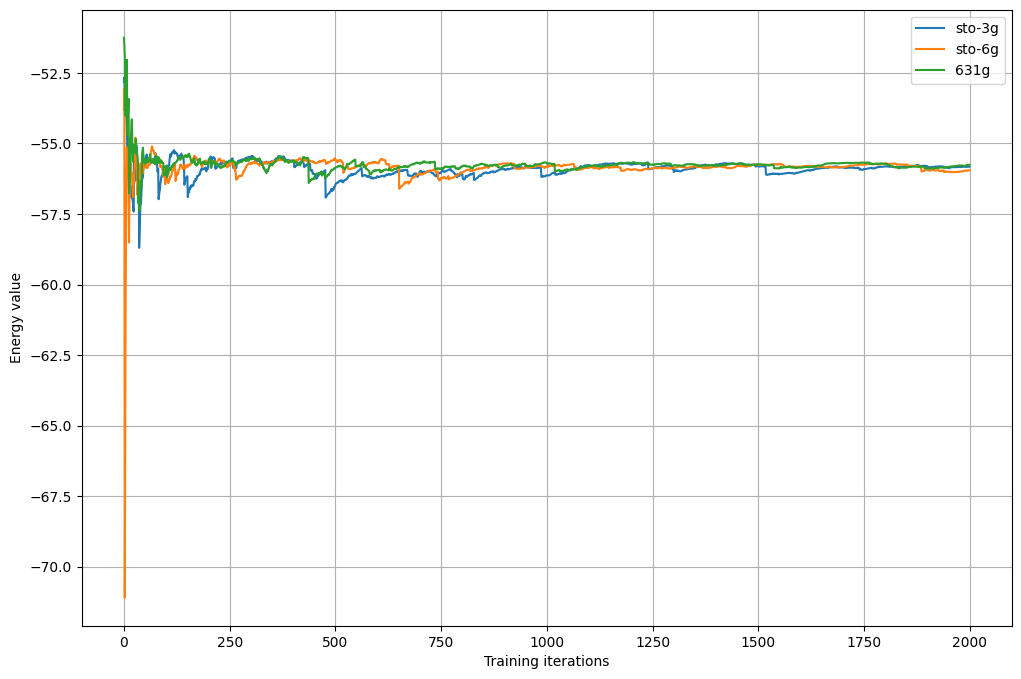}
    \caption{Paulinet : energy minimization of $NH_3$ with different basis sets}
    \label{fig:pretraining_pn_nh3}
\endminipage
\end{figure}

PauliNet on the other side makes use of the single electron orbitals, which requires a choice of basis set. We investigate how different basis sets impact the convergence (Figure~\ref{fig:pretraining_pn_c2h4} and Figure~\ref{fig:pretraining_pn_nh3}).

\begin{figure}[H]
\minipage{0.5\textwidth}
    \includegraphics[width=\textwidth]{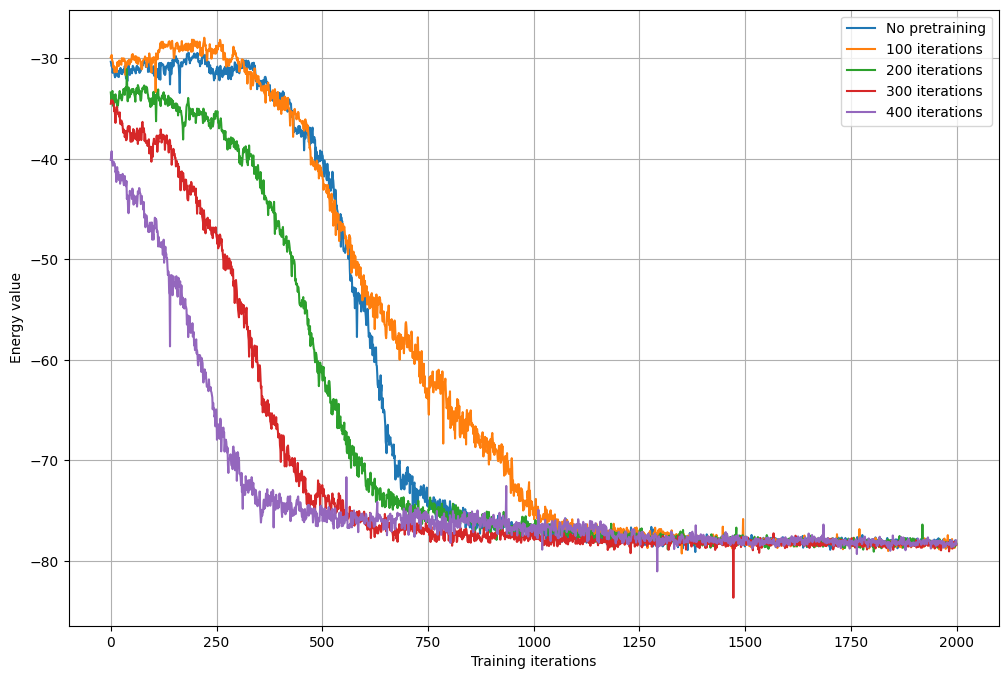}
        \caption{Ferminet : energy minimization of $C_2H_4$ for different pretraining iterations}
        \label{fig:pretraining_fn_c2h4}
\endminipage
\minipage{0.5\textwidth}
\includegraphics[width=\textwidth]{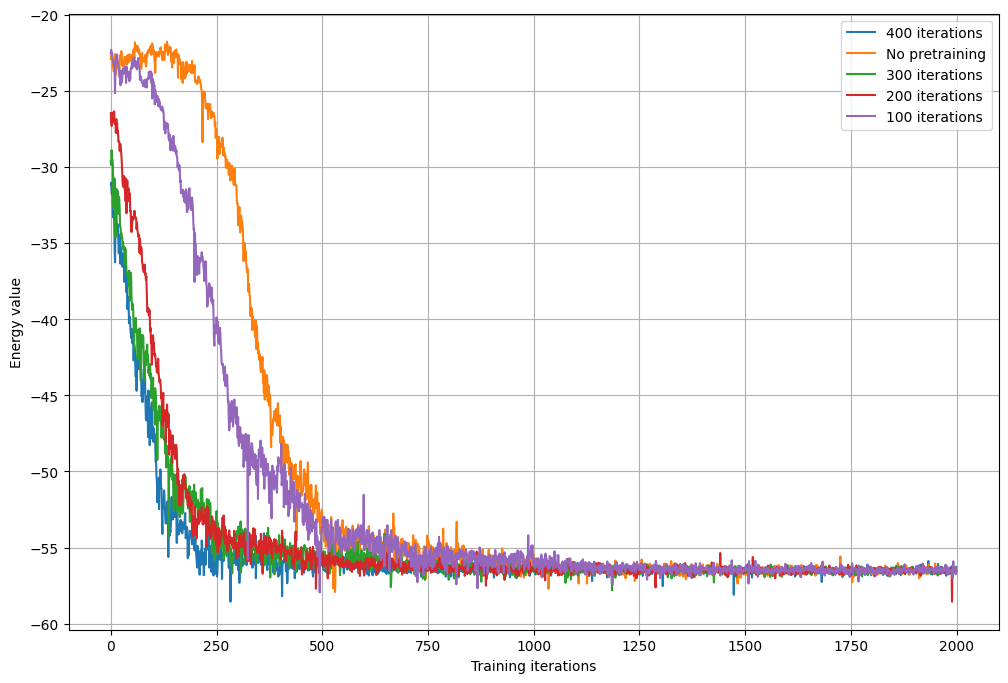}
    \caption{ Ferminet : energy minimization of $NH_3$ for different pretraining iterations}
    \label{fig:pretraining_fn_nh3}
\endminipage
\end{figure}

We observe (Figure~\ref{fig:pretraining_fn_c2h4} and Figure~\ref{fig:pretraining_fn_nh3})that the pretraining step in FermiNet is of tremendous importance to accelerate the learning phase.  An interesting observation is the fact that pretraining the model for $ C_2H_4 $ with only $ 100 $ iterations leads to slower convergence towards the ground state than no pretraining, which is probably due to the fact that the model has not been able to reach the desired target (eigenstate given by Hartree-Fock) in 100 iterations.
Nonetheless, we point out that the learned state remains the same after a sufficiently large number of iterations.

\subsubsection{Sampling points}

FermiNet and PauliNet are trained using variational Monte Carlo to evaluate high-dimensional integrals. Training data are generated from $ p(\bfm{r}) \sim \psi^2_{\theta}(\bfm{r})$ after a certain number of training iterations using the latest $ \psi_{\theta}(\bfm{r})$. A Monte Carlo Markov Chain (MCMC) algorithm is required in order to sample electrons configurations from probability distribution $ p $ (Metropolis-Hastings for ferminet and Langevin Monte Carlo for Paulinet) . We investigate how the number of sampling points impact the training.

\begin{figure}[H]
\minipage{0.5\textwidth}
    \includegraphics[width=\textwidth]{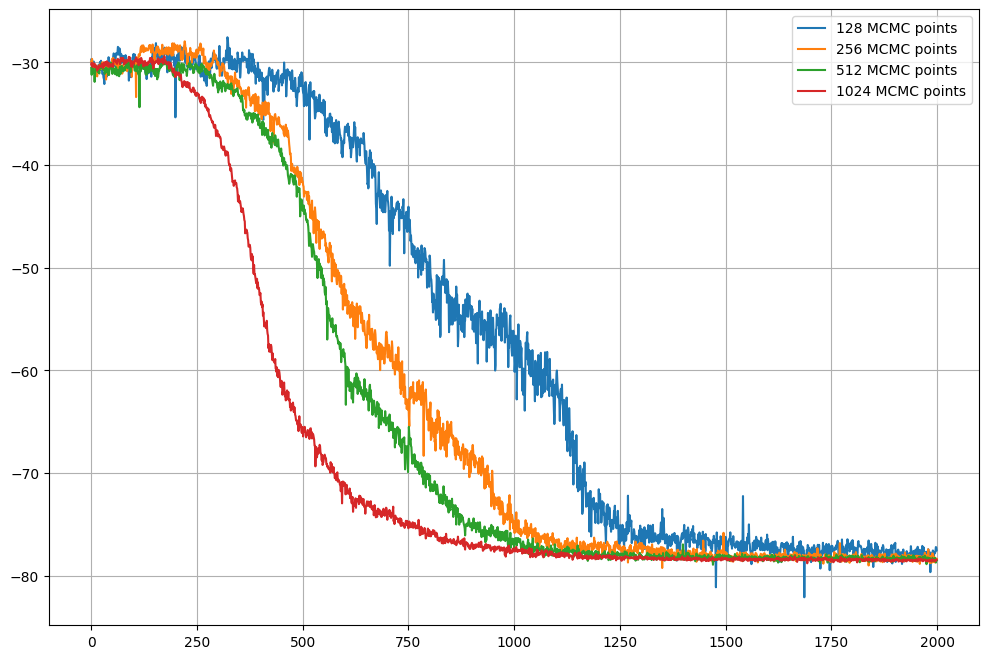}
        \caption{ Ferminet : energy minimization of $C_2H_4$ for different number of integration points}
        \label{fig:ferminet_c2h4_mcmc}
\endminipage
\minipage{0.5\textwidth}
\includegraphics[width=\textwidth]{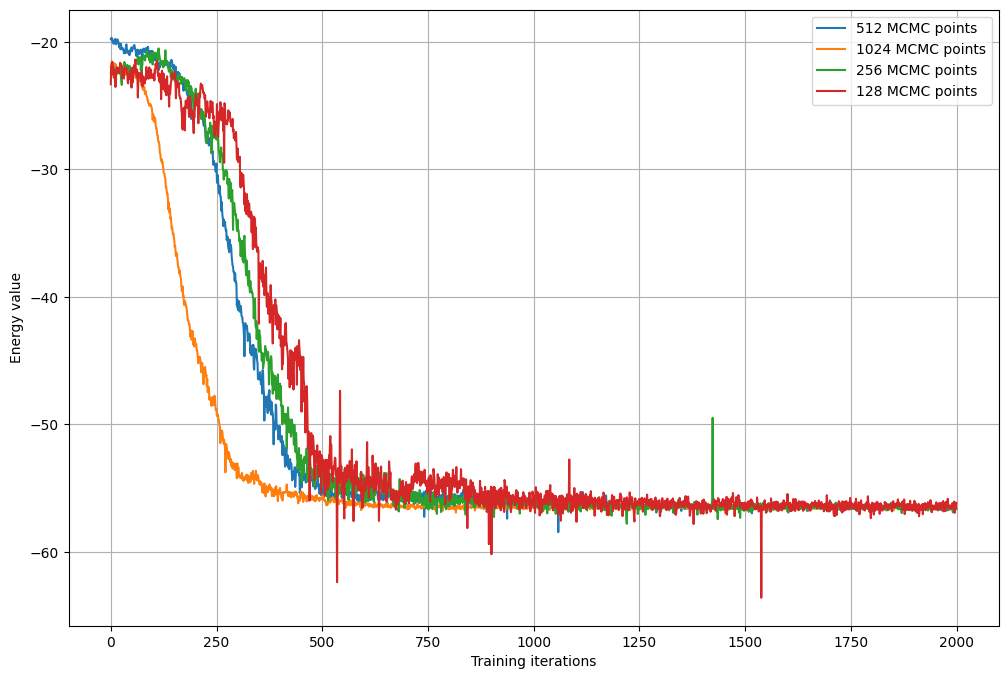}
    \caption{Ferminet : energy minimization of $NH_3$ for different number of integration points}
    \label{fig:ferminet_nh3_mcmc}
\endminipage
\end{figure}

\begin{figure}[H]
\minipage{0.5\textwidth}
    \includegraphics[width=\textwidth]{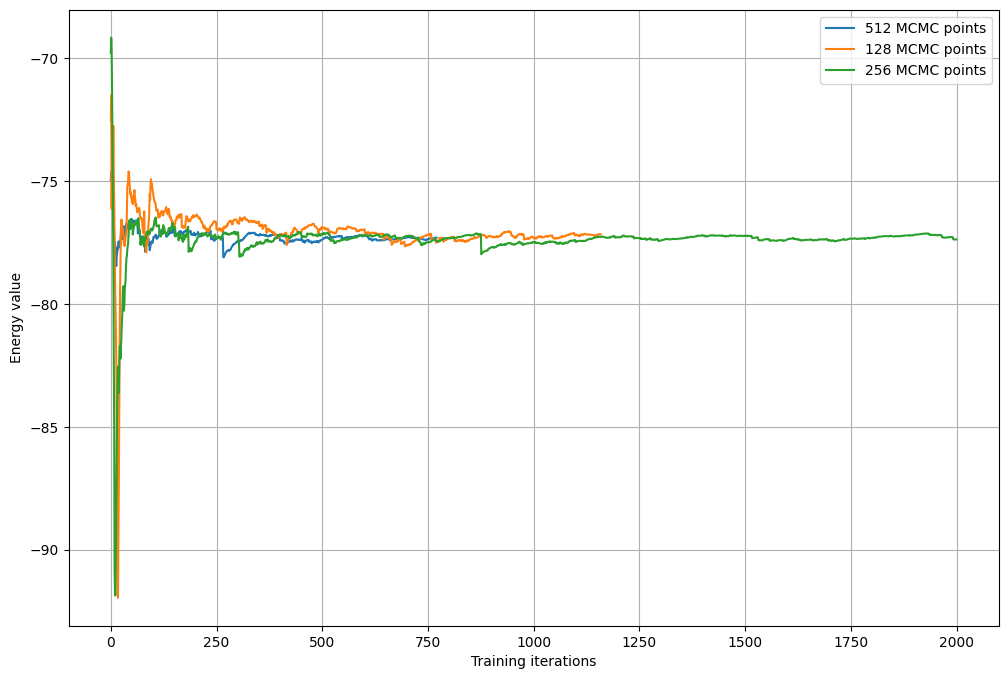}
        \caption{ Paulinet : energy minimization of $C_2H_4$ for different number of integration points}
        \label{fig:paulinet_c2h4_mcmc}
\endminipage
\minipage{0.5\textwidth}
\includegraphics[width=\textwidth]{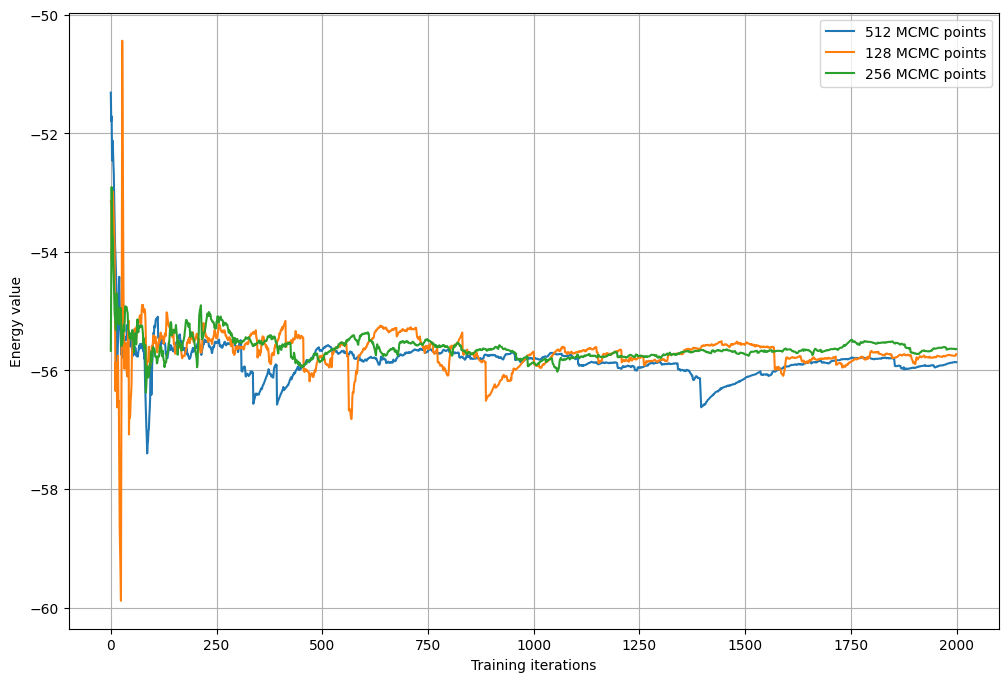}
    \caption{Paulinet : energy minimization of $NH_3$ for different number of integration points}
    \label{fig:paulinet_nh3_mcmc}
\endminipage
\end{figure}

The results (Figures \ref{fig:ferminet_c2h4_mcmc}, \ref{fig:ferminet_nh3_mcmc}, \ref{fig:paulinet_c2h4_mcmc},\ref{fig:ferminet_nh3_mcmc}) show that the more samples we have, the faster the training for Ferminet. The final learned state though, does not depend on our choice on the number of samples.
\subsubsection{Number of determinants}
\begin{figure}[H]
\minipage{0.5\textwidth}
    \includegraphics[width=\textwidth]{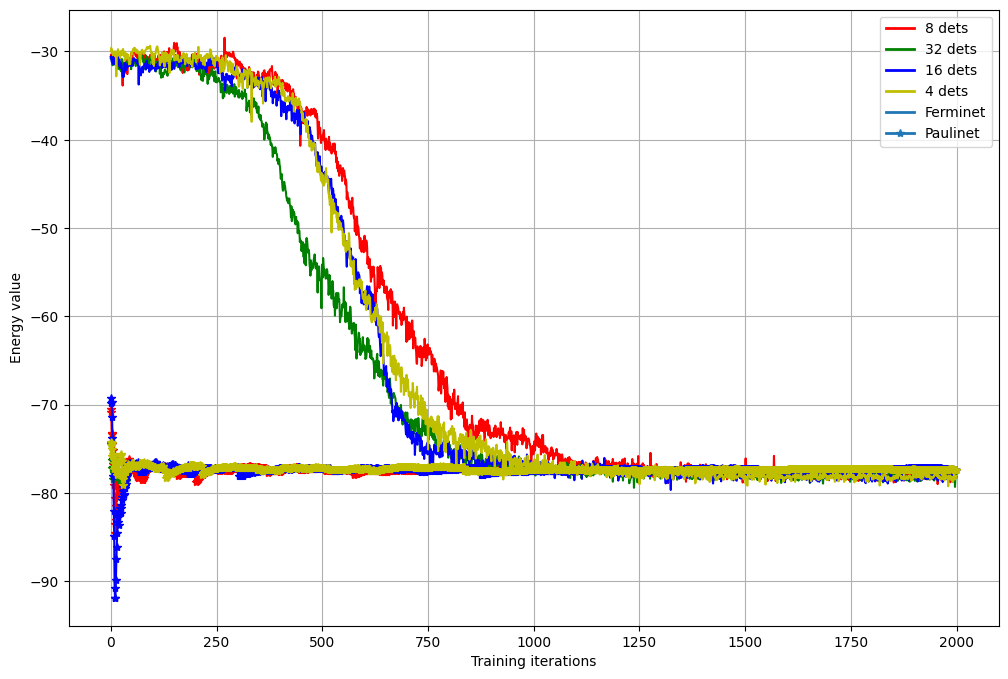}
        \caption{ Energy minimization of $C_2H_4$ with different number of dets}
        \label{fig:nb_det_c2h4}
\endminipage
\minipage{0.5\textwidth}
\includegraphics[width=\textwidth]{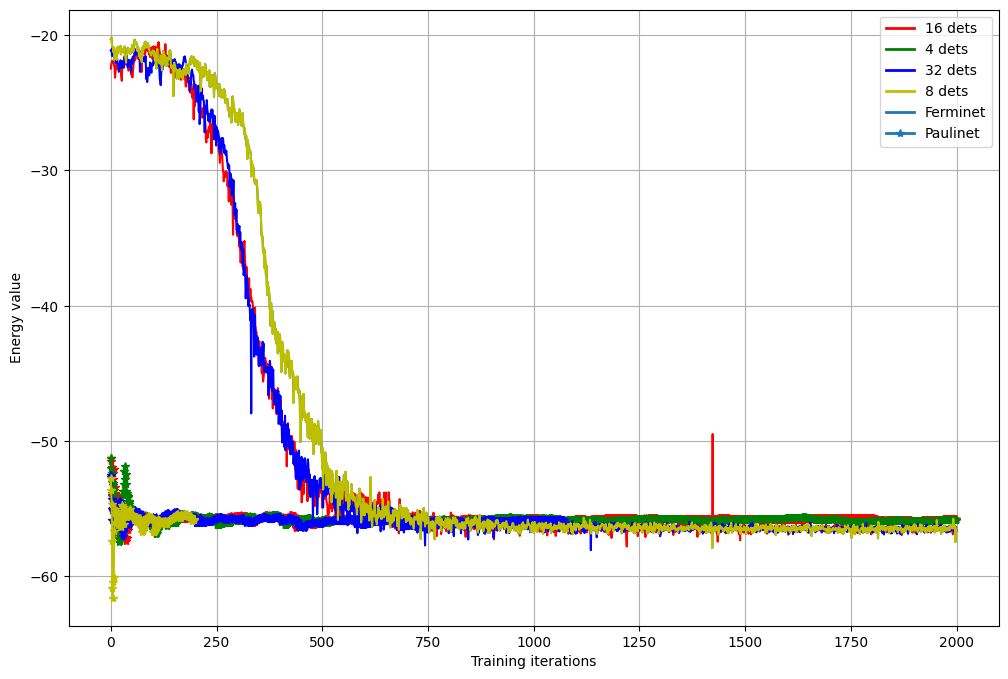}
    \caption{Energy minimization of $NH_3$ with different number of dets}
    \label{fig:nb_det_nh3}
\endminipage
\end{figure}

Quantum chemistry algorithms generally require many Slater determinants to reach high accuracy. The results here (Figures \ref{fig:nb_det_c2h4}, \ref{fig:nb_det_nh3}) show how neural networks can help to reduce the number of determinants needed. We observe convergence for both models on both test cases, even with only 4 determinants.   
\subsubsection{Optimizers for Ferminet}
We finally illustrate the importance in choice of the gradient descent algorithm for ferminet. We compare Adam optimizer with KFAC, which is a second order optimization algorithm derived from natural gradient descent. KFAC outperforms Adam for both molecules. 
\begin{figure}[H]
\minipage{0.5\textwidth}
    \includegraphics[width=\textwidth]{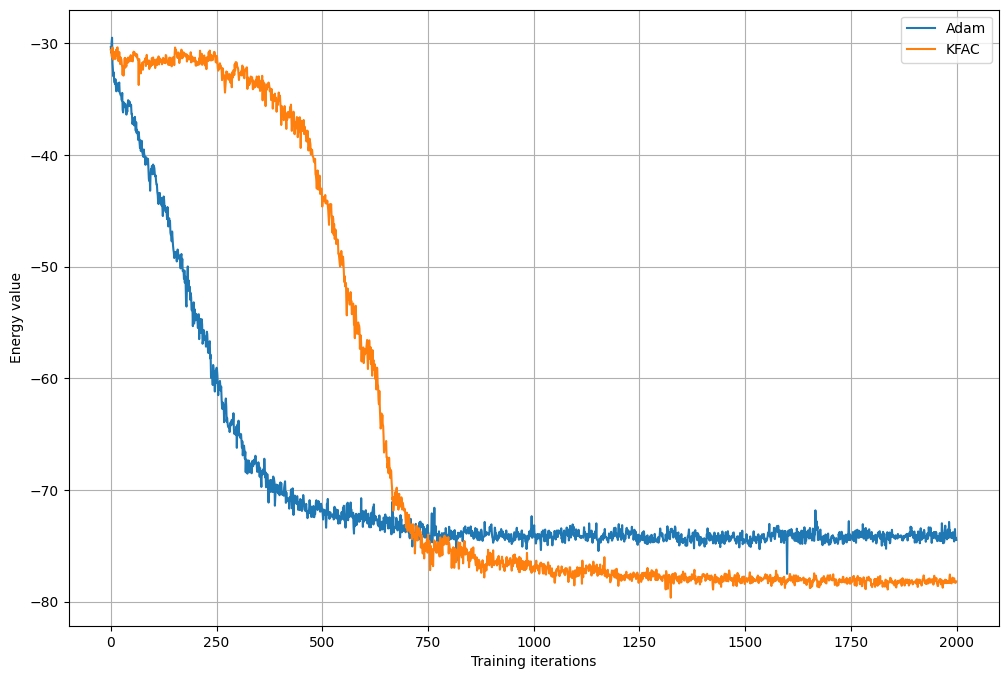}
        \caption{ Ferminet : Energy minimization of $C_2H_4$ with Adam and KFAC}
        \label{fig:vqmc_sampling}
\endminipage
\minipage{0.5\textwidth}
\includegraphics[width=\textwidth]{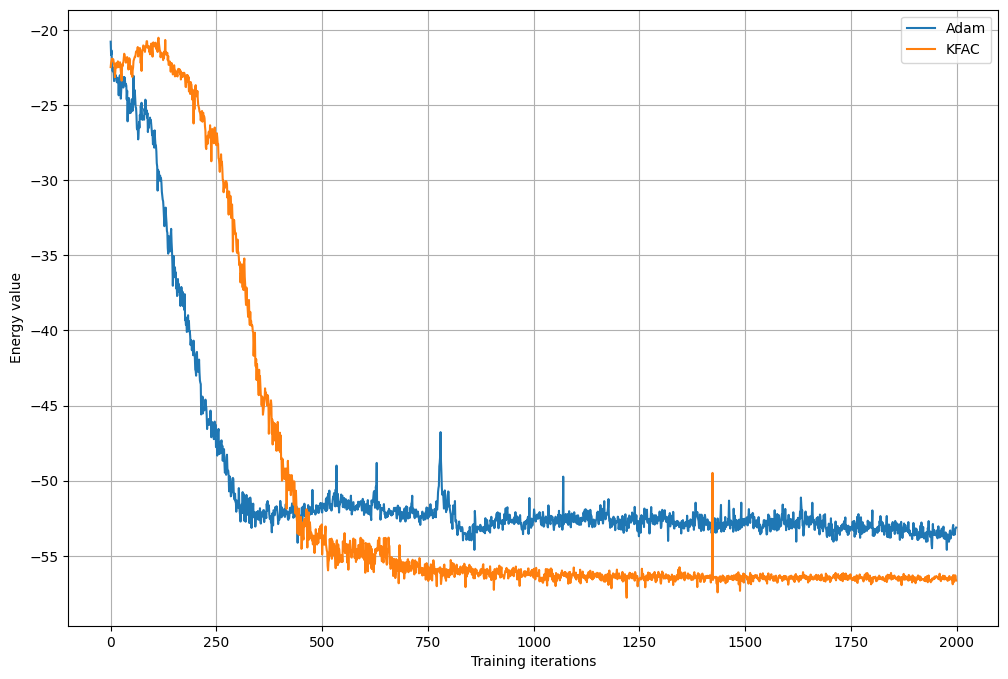}
    \caption{Ferminet : Energy minimization of $NH_3$ with Adam and KFAC}
    \label{fig:mean_variance_pn}
\endminipage
\end{figure}
\subsection{Comparison between the different methods presented above} \label{sec: Numerical experiments}

We display on Table~\ref{tab:tab1} the different values of the energy obtained for the different methods with several molecules. The exact values were obtained from \cite{PhysRevResearch.2.033429}.
FermiNet was trained using a $ L = 4$ layers neural networks, and activation function $ \sigma = tanh $. We used $ 256 $ training samples at each iteration, and summed over $ P = 16 $ determinants. PauliNet on the other side, was modeled with a  3-neural networks layers for the Jastrow factor and backflow. There were 256 training samples used at each iteration, and the wavefunction was approximated summing over 16 determinants.

Our simulations suggest that FermiNet gives in general better approximation results than PauliNet. Indeed, the energy prediction of the $ LiH $ and the oxygen by FermiNet is considerably better. This however comes at a cost, as the number of parameters optimized in FermiNet is higher than that of PauliNet.
There is therefore a tradeoff between accuracy and complexity between these two methods, which needs to be investigated also for larger molecules.

For DMRG, the rank was taken equal to $6$, increasing the rank further did not lead to any improvement on these examples. The molecule $C_2H_4$ could not be computed, since even the minimal basis would have required at least $16$ functions, which was too expensive to be computed on the computers available to us at that moment. The high cost of the computation is due to the construction of the TTO representing the Hamiltonian operator $\tilde H$ defined by \eqref{equ: Fock hamiltonian}. Indeed, to avoid having to deal with very large cores in the TTO, the algorithm uses a "compression algorithm" which leads to an approximation of a given TTO with smaller cores, however, this compression ends up being quite expensive as the cores grow, which is why we were so limited for the computation.

Compared to PauliNet and FermiNet, DMRG seems to give less accurate results. But this probably comes from the use of very small basis sets.

\begin{table}[h!]
    \centering
\begin{tabular}{cccc|cccc|cccc|cccc|cccc}
\hline
    Atom / Molecule &  DMRG (Energy, Basis) & FermiNet (Energy) & PauliNet (Energy) & Exact\\
\hline
    Li & -7.431553 (6-31G) & -7.473162 & \textbf{-7.4762} & -7.47806032\\
\hline
    C & -37.673536 (6-31G) & -37.8077 & \textbf{-37.8176} & -37.8450\\
\hline
    O & -74.516816 (STO-6G) & \textbf{-75.0551} & -75.008 &  -75.0673\\
\hline
    $H_2$ & -1.014310 (6-31G) & -1.1378366 & -1.02192 & $\boldsymbol \times$ \\
\hline
    LiH & -7.949315 (6-31G) & \textbf{-8.073885} & -8.00654 & -8.070548\\
\hline
    $NH_3 $ & -55.713815 (STO-6G) & -56.1439 & -56.7146 & $\boldsymbol \times$\\
\hline
    $C_2H_4$ & $\boldsymbol \times$  & -78.52951 & -76.0121 & $\boldsymbol \times$\\
\hline
\end{tabular}
    \caption{Comparison of the energy computed for the different methods for several atoms and molecules}
    \label{tab:tab1}
\end{table} 

\begin{table}[h!]
    \centering
\begin{tabular}{cccc|cccc|cccc|cccc}
\hline
    Atom / Molecule &  DMRG (time (s), Basis) & FermiNet (time (s)) & PauliNet (time (s))\\
\hline
    Li & 20.3 (6-31G) & 47.2 & 108.1\\
\hline
    C & 31.8 (6-31G) & 57.5 & 137.9\\
\hline
    O & 1.1 (STO-6G) & 70.1 & 166.5\\
\hline
    $H_2$ & 0.75 (6-31G) & 37.4 & 72.7 \\
\hline
    LiH & 87.6 (6-31G) & 45 & 117.6 \\
\hline
    $NH_3 $ & 16.4 (STO-6G) & 71.3 & 260.3 \\
\hline
    $C_2H_4$ & $\boldsymbol \times$ & 146.9 & 590.6 \\
\hline

\end{tabular}
    \caption{Comparison of the computation time of the different methods for several atoms and molecules}
    \label{tab:tab2}
\end{table}

\section{SimplexNet}
\label{sec:simplexnet}

Previous methods and additional ones in the litterature \cite{Lin2023} are based on using determinants to ensure the antisymmetry condition of the wavefunction $\psi$. 
This can be heavy theoretically with the second quantization and even computationally for the training of neural networks presented in previous sections. 
Therefore, we now present an alternative approach which avoids the use of these determinants. 

For pedagogical reasons, we present the method on a simplified problem. Namely, we restrict ourselves to the domain $[0,1]^{N \times \kappa}$ and do not consider spin related effects, which is not physical but easier to analyze theoretically. Here the variable $\kappa$ corresponds the dimension in which particles moves and in particular, it could be different from three.

Let $S_N$ be the group of permutations with $N$ elements,
and define $S := \{ (\bfm{r}_i)_{1\le i \leq N} \in [0,1]^{N \times \kappa} \ | \ | \bfm{r}_1| \leq \cdots \leq | \bfm{r}_N| \}$ be what we call the ``principal simplex''. One can easily define the projection $\Pi$ over $S$, the set $S$ being convex and closed. As the wavefunction $\psi$ is antisymmetric, there holds

\begin{equation}\label{eq:anzatz}
\psi = \varepsilon(\Pi) \psi \circ \Pi.
\end{equation}

and consequently, it is sufficient to know the value of $\psi$ on the principal simplex $S$ to evaluate $\psi$ everywhere. However contrary to the initial multi-body Schr\"odinger problem posed on $\R^{3N}$, here we have a boundary condition on the boundary of $S$ which is of Dirichlet type that is

$$
\psi = 0 \text{ on } \partial S.
$$
This is automatically satisfied with our anzatz \ref{eq:anzatz}.

\begin{figure}[H]
    \centering
    \begin{subfigure}[b]{0.4\textwidth}
        \centering
        \includegraphics[scale=0.5]{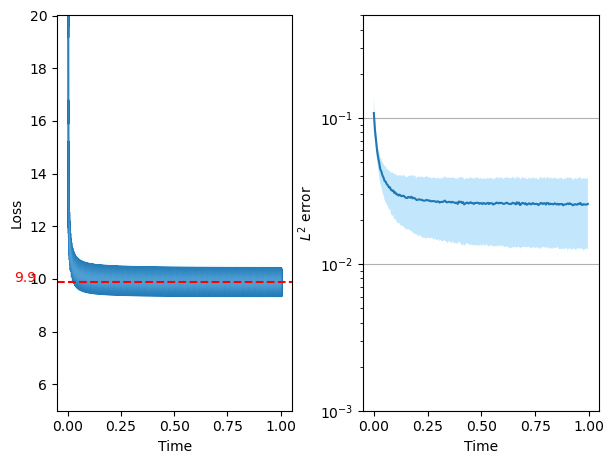}
        \caption{Plot of the energy and the $L^2$ error}
        \label{fig:energySimplesNet}
    \end{subfigure}
    \hfill
    \begin{subfigure}[b]{0.4\textwidth}
        \centering
        \includegraphics[scale=0.5]{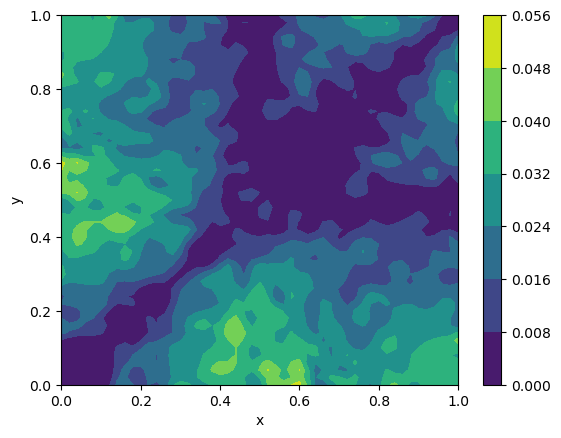}
        \caption{Plot of the local error}
        \label{fig:solSimplesNet}
    \end{subfigure}
    \caption{Test of SimplexNet on a simple case}
    \label{fig:resultsSimplexNet}
\end{figure}

We use the idea in the context of neural networks adding a layer we call the projection layer at the beginning of a neural network representing $\psi$. In this section, we will only use feedforward neural networks and the global architecture called ``SimplexNet'' is represented on Figure~\ref{fig:energySimplesNet} where the network is composed by the input layer, the projection layer, a stack of dense layers (in Figure ~\ref{fig:energySimplesNet} there is only one hidden layer but one can add many more) and the output layer. Note that in order to preserve the antisymmetry of the represented function, the signature of the permutation $s$ such that

\[
(\mathbf{r}_{s(1)}, \cdots, \mathbf{r}_{s(N)}) := \Pi (\mathbf{r}_1, \cdots, \mathbf{r}_N).
\]
is used to compute the output. Additionally and taking the example of FermiNet \cite{Pfau_2020}, a forward/backward pass with a layer computing a single determinant is of complexity $N^3$ as an singular value decomposition is calculated. For our architecture, a forward pass costs $O(N \log(N))$ operations since a sorting routine is used. For the backward pass, things are a bit more complicated since the sorting operation is not everywhere differentiable. To deal with this difficulty, we forced the jacobian matrix of this layer to be equal to the matrix $M_s$ defined by:

$$
M_{s,i,j} = 
\left\{
\begin{array}{rl}
  1  & \text{ if j = s(i)} \\
  0  & \text{ otherwise.}
\end{array}
\right.
$$
which coincides with the jacobian of the sorting operator when it is smooth.
Hence the complexity of a backward pass is $O(N)$. The signature is approximated by the differentiable function :

$$
\varepsilon(s) := \prod_{1 \leq i < j \leq N } \tanh(|\bfm{r}_j| - |\bfm{r}_i|)
$$
which takes $O(N^2)$ operations.

\tikzset{%
  every neuron/.style={
    circle,
    draw,
    minimum size=0.1cm
  },
  neuron missing/.style={
    draw=none, 
    scale=1,
    text height=0.05cm,
    execute at begin node=\color{black}$\vdots$
  },
}

\begin{figure}[h]
\centering
\begin{tikzpicture}[scale=0.95, >=stealth]

\foreach \m/\l [count=\y] in {1,2,3,missing,4}
  \node [every neuron/.try, neuron \m/.try] (input-\m) at (0,2.5-\y) {};

\foreach \m/\l [count=\y] in {1,2,3,missing,4}
  \node [every neuron/.try, neuron \m/.try] (project-\m) at (2,2.5-\y) {};  

\foreach \m [count=\y] in {1,missing,2}
  \node [every neuron/.try, neuron \m/.try ] (hidden-\m) at (4,2-\y*1.25) {};

\foreach \m [count=\y] in {1}
  \node [every neuron/.try, neuron \m/.try ] (output-\m) at (6,0.75-\y) {};

\foreach \l [count=\i] in {1,2,3,N}
  \draw [<-] (input-\i) -- ++(-1,0)
    node [above, midway] {$\mathbf r_\l$};

\foreach \l [count=\i] in {1,2,3,N}
  \node [above] at (project-\i.north) {$\mathbf r_{s(\l)}$};
  
\foreach \l [count=\i] in {1,m}
  \node [above] at (hidden-\i.north) {$H_\l$};

\foreach \l [count=\i] in {1}
  \draw [->] (output-\i) -- ++(1,0)
    node [above=0.3cm] {$\psi(\mathbf{r}) = \varepsilon(s) \psi(\Pi \mathbf{r})$};
    
\foreach \i in {1,...,4}
  \foreach \j in {1,...,4}
    \draw [->] (input-\i) -- (project-\j);
    
\foreach \i in {1,...,4}
  \foreach \j in {1,...,2}
    \draw [->] (project-\i) -- (hidden-\j);

\foreach \i in {1,...,2}
  \foreach \j in {1}
    \draw [->] (hidden-\i) -- (output-\j);

\foreach \l [count=\x from 0] in {Input, Projection, Hidden, Output}
  \node [align=center, above] at (\x*2,2) {\l \\ layer};

\end{tikzpicture}
\caption{SimplexNet}
\label{fig:SimplexNet}
\end{figure}
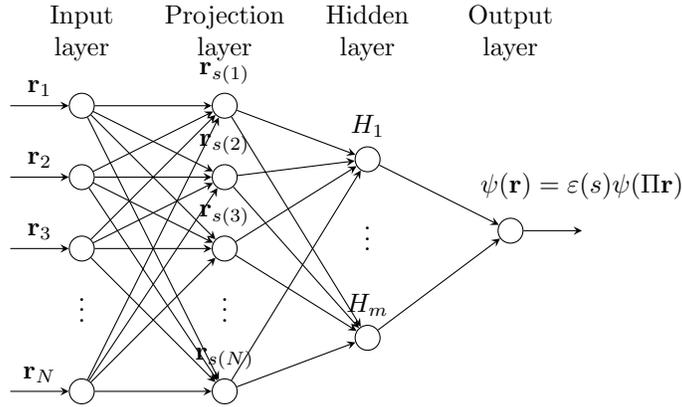

In order to provide a first test of the SimplexNet, we will not consider the full Schr\"odinger multi-body problem as it involves an unbounded domain and high dimensionality. To remove these difficulties,  we consider instead the Poisson--Neumann problem on a bounded domain $\Omega := [0,1]^N$. In such case, the hamiltonian is just $-\frac 1 2 \sum_{j=1}^N \Delta_{\bfm{r}_j}$.

We use the  framework \href{https://www.tensorflow.org/?hl=fr}{python/Tensorflow/Keras} such that :

\begin{itemize}
    \item The network has one hidden layers where the first one has $10^3$ neurons.
    \item We use the stochastic gradient descent as optimizer so that the minimized energy writes

    $$
    \frac{1}{n} \sum_{i=1}^n \sum_{j=1}^N |\nabla_{\mathbf r_j} \psi(\mathbf{r}_i)|^2 
    $$
    where the $(\mathbf{r}_i)$ are sampled in the domain $\Omega$. Note that this implies that we have one dataset and the batch size is $n=1000$. 
    
    \item The output layer is linear and normalized at each iteration step to have a unit $L^2$ norm solution.
\end{itemize}

In Figure \ref{fig:resultsSimplexNet},  we plot our results for $N=2$ so that the exact solution is $\psi(\bfm{r}_1,\bfm{r}_2) = \pm (\cos(\pi \bfm{r}_1) - \cos(\pi \bfm{r}_2))$ and the fundamental energy is $\pi^2$.

The time in Figure \ref{fig:energySimplesNet} corresponds to the number of iterations multiplied by the time step used in the stochastic gradient algorithm. In the same picture, the training was run 5 times such that the thick line corresponds to mean values while the shaded area corresponds to variance. It can be noticed that such variance is relatively small at the end of the training process.

To summarise, the method is indeed efficient for such simple problem. For more physical relevant problems, the additional difficulties of an unbounded domain and a singular coulombic potential prevents us to get relevant results.
One should add more structure to the neural network exploiting and adapting the anzatz of FermiNet and PauliNet.

\section{Conclusion}
\label{sec:concl}

In this paper, we presented a comparison of tensors and neural network methods to solve the multibody Schr\"odinger eigenvalue problem. Neural network methods seems to be at least as precise as DMRG. These methods are mainly empirical and no proof of convergence of those algorithms seems to be available. In addition, another neural network approach called SimplexNet is explored preventing the use of expensive determinants. It was tested successfully on a simple non singular case and has a theoretical justification contrary to other methods \cite{DusEhrlacher2024}. However for the moment the neural network involved is over-simplistic and does not encode enough physical effects to solve the very singular Schr\"odinger problem. Hence, it should be interesting to see if ideas taken from FermiNet and PauliNet can be included into the SimplexNet method.

\bibliographystyle{siam}
\bibliography{refs}

\end{document}